\newcommand{\argmin}{\operatornamewithlimits{argmin}}

\def\qed{\hfil\hfil{\vrule height7pt width6pt depth0pt}\hfil}

\documentclass[12pt]{article}
\newtheorem{theorem}{Theorem}[section]

\newtheorem{lemma}[theorem]{Lemma}

\newtheorem{remark}[theorem]{Remark}

\usepackage{amssymb}
\usepackage{amsmath}
\usepackage{graphicx}% Include figure files in eps format.
\usepackage{bm}% bold math.
\usepackage{epstopdf}
\usepackage[ruled,vlined,linesnumbered]{algorithm2e}
\usepackage{enumerate}
\usepackage{multirow}

\def\proof{\noindent{\bf Proof: }}
\def\arg{{\hbox{arg}}}

\def\min{{\hbox{min}}}
\def\sgn{{\hbox{sgn}}}
\def\circulant{{\hbox{circ}}}

\def\B{{\mathcal{B}}}
\def\d{~\mathrm{d}}
\def\x{\mathbf{x}}

\def\w{\mathbf{w}}
\def\j{\mathbf{j}}
\def\dx{~\mathrm{d}\mathbf{x}}

\def\SOPW{\mathbf{SOPW}}
\def\ISOPW{\mathbf{ISOPW}}
\def\SOBF{\mathbf{SOBF}}
\def\ISOBF{\mathbf{ISOBF}}

\title{Projection to the Set of Shift Orthogonal Functions\thanks{The research of R.C. and F.B. was supported by the Department of Energy grant number
DOE-SC0010613. The research of S.O. and K.Y. was supported by the Office of Naval Research (Grant N00014-11-1-719). The research of V.O. was supported by the National Science Foundation under Award DMR-1106024}}

\author{
Farzin Barekat\thanks{Department of Mathematics, University of California, Los Angeles, USA.} ~
,
Rongjie Lai \thanks{Department of Mathematics, University of California, Irvine, USA.}~
,
Ke Yin  \thanks{Department of Mathematics, University of California, Los Angeles, USA.} ~
,
Stanley Osher\thanks{Department of Mathematics, University of California, Los Angeles, USA.} 
,
\\
Russel Caflisch\thanks{Department of Mathematics, University of California, Los Angeles, USA.} 
,
Vidvuds Ozoli{\c{n}}{\v{s}}\thanks{Department of Materials Science and Engineering, University of California, Los Angeles, USA.}
}
\begin{document}
\maketitle

\begin{abstract}
This paper presents a fast algorithm for projecting a given function to the set of shift orthogonal functions (i.e. set containing functions with unit $L^2$ norm that are orthogonal to their prescribed shifts). The algorithm can be parallelized easily and its computational complexity is bounded by $O(M\log(M))$, where $M$ is the number of coefficients used for storing the input. To derive the algorithm, a particular class of basis called Shift Orthogonal Basis Functions are introduced and some theory regarding them is developed.  
\end{abstract}

\section{Introduction}\label{section:introduction}

%In this manuscript we introduce a class of basis, which we call Shift Orthogonal Plane Waves (SOPWs), that have some desirable properties both for theory and in practice. In particular we advocate that in the variational problems with shift orthogonality constraints, using these basis improves numerical calculations both in speed and accuracy. 

Let  $\Omega=[0,L_1]\times \cdots\times [0,L_d]$ be a bounded, periodic domain.  Let $\w = (w_1,\cdots,w_d)\in\mathbb{R}_{+}^d$ be a basis of a $d$-dimensional lattice, and define
\begin{equation}
\label{eqn:lattice}
\Gamma_{\w} = \{\j\w:=(j_1w_1,\cdots,j_dw_d)\hspace{0.1cm} |\hspace{0.1cm}\j=(j_1,\cdots,j_d)\in\mathbb{Z}^d\}.
\end{equation}
Note that to make the boundary of $\Omega$ periodic, we require that $L_i$ be divisible by $w_i$, for $i=1,\ldots,d$.  Endow $\Omega$ with the usual inner product
\[ \langle f,g\rangle=\int_\Omega f^*g \dx.\]
For any function $f(\x)\in L^2(\Omega)$, we say $f$ is \emph{shift orthogonal}, if it satisfies \emph{shift orthogonality constraints}:
\begin{equation} \int_\Omega f(\x)^*f(\x-\j\w) d\x=\delta_{\j0} \qquad \hbox{ for all }  \j\w\in \Gamma_{\w}. \label{eqaution:shiftOrthogonalDefinition}\end{equation} 
Shift orthogonality constraints arise naturally in many applications of science and engineering. However, due to the numerical and theoretical challenges in imposing the shift orthogonality constraints, these constraints have not received much attention in the literature. 

%The main contribution of this paper is to propose a very fast algorithm to find projection of a given functions to the set of shift orthogonal functions (i.e. to find a shift orthogonal functions that is closest in $L^2$ norm to the given functions). 

The main contribution of this paper is to propose a very fast algorithm for finding a shift orthogonal functions that is closest in $L^2$ norm to a given function (i.e. projection to the set of shift orthogonal functions). Note that shift orthogonal functions constitute a set and not a vector space. There are many potential applications for the algorithm described in this paper. As an example, in section \ref{section:applicationCPW} we demonstrate how the algorithm increases the speed in computing CPWs described in \cite{Ozolins:2013CPWs}. Further applications of the algorithm will be investigated in subsequent research projects.

In order to devise the algorithm we introduce the notion of Shift Orthogonal Basis Functions (SOBFs) and develop some theory for them. For one dimensional periodic domain $[0,L_1]$ we call a family of functions
\[ \{\xi^i_j(x)\}_{i=1,j=0}^{i=\infty,j=L_1/w_1-1},  \]
with $\xi^i_j(x)=\xi^i_0(x-jw_1)$, \emph{Shift Orthogonal Basis Functions (SOBFs)} if they form a complete orthonormal set of basis for space $L^2([0,L_1])$. We refer to superscript index $i$ as the \emph{depth index}, and subscript index $j$ as the \emph{shift index}. For higher dimensional periodic domains, we call a complete set of basis SOBFs if they are formed from the tensor product of one dimensional SOBFs. Section \ref{section:definitionSOPW} describes a concrete example of SOBFs with certain nice properties.

It turns out that SOBFs have been of interest to quantum mechanics and signal analysis communities (i.e. although, the orthonormal basis considered in these literatures are usually defined in domain $\mathbb{R}^d$, instead of periodic bounded domains).  In \cite{BIB:Wilson,BIB:Sullivan}, a numerical example of SOBFs, called phase space Wannier functions, are provided that are exponentially localized in time and frequency domain and for which the matrix of the Laplacian operator is near-diagonal (i.e. see \cite[equations 2.1, 2.3 and 2.5]{BIB:Sullivan}). Reference \cite{BIB:Daubechies} provides a simple description of certain SOBFs, called Wilson basis, whose Fourier transform have specific bimodal form (see \cite[equations 1.9a and 1.9b]{BIB:Daubechies}).  An explicit example of Wilson basis with exponential decay in time and frequency domain is also constructed in \cite{BIB:Daubechies} using rapidly converging superpositions of Gaussians  (i.e. however, the Laplacian matrix is no longer near-diagonal for this basis). Other examples of Wilson basis are presented in \cite{BIB:Laeng} and \cite{BIB:Auscher}.  

%It turns out that some specific type of SOBFs, called Wilson basis \cite{BIB:Wilson,BIB:Sullivan,BIB:Daubechies}, have been of interest to quantum mechanics and signal analysis communities. Although, Wilson basis are defined for domain $\mathbb{R}^d$, the adapted version of them for periodic bounded domains belongs to class SOBFs. Indeed, Wilson basis are SOBFs whose Fourier transform have specific bimodal form (see \cite[equations 1.9a and 1.9b]{BIB:Daubechies}). A numerical example of Wilson basis, called Phase plane Wannier Functions, were constructed in \cite{BIB:Sullivan} that were exponentially localized both in time and frequency domain and the Laplacian matrix is near-diagonal. In  \cite{BIB:Daubechies} also provide an explicit example of Wilson basis using rapidly converging superpositions of Gaussians that are easy to construct and have exponential decay in time and frequency domain (i.e. however, the Laplacian matrix is no longer near-diagonalized for these basis). 

The remainder of this paper consists of the following: Section \ref{section:theorySOBF2D} develops some theory regarding SOBFs that is used in devising the fast algorithm. Section \ref{section:projectionAlgorithm} contains the description of the fast algorithm and highlights some of its important properties including its computational complexity.  Section \ref{section:definitionSOPW} describes a concrete example of SOBFs, which we call Shift Orthogonal Plane Waves (SOPWs), that have certain nice properties. Section \ref{section:applicationCPW} gives an application of the fast algorithm in generating CPWs. Section \ref{section:conclusion} presents some concluding remarks. Appendices \ref{section:variation}, \ref{section:sameDepth} and \ref{section:firstSecondDerivative} contain proofs for properties of SOPWs mentioned in section \ref{section:definitionSOPW}.

\section{Shift Orthogonal Basis Functions (SOBFs)}\label{section:theorySOBF2D}
This section develops some theory for SOBFs in multidimensional domains. Some of this theory also appears in \cite{BIB:BarekatThesis}. For the ease of exposition, the concepts are illustrated for 2D domain. The results to other dimensions can be extended in the trivial way. 

By scaling, it is  assumed without loss of generality that $\Omega=[0,L_1]\times [0,L_2]$ and the length of the shift in each coordinate is unit length (i.e. $L_1$ and $L_2$ are replaced by $L_1/w_1$ and $L_2/w_2$, respectively). Because SOBFs in 1D form an orthonormal set of basis, tensor product can be used to form orthonormal set of basis for 2D. That is, any function $f\in L^2(\Omega)$ can be expressed by 
\begin{equation} f(x_1,x_2)=\sum_{i_1,i_2=1}^{\infty}\sum_{j_1=0,j_2=0}^{L_1-1,L_2-1} a_{j_1,j_2}^{i_1,i_2}\xi_{j_1}^{i_1}(x_1)\xi_{j_2}^{i_2}(x_2).\label{equation:fSOBFexpansion2D}\end{equation}

\begin{remark}  
In what follows $i_1$ and $i_2$ and their primes (i.e. $i'_1$, $i''_2$, etc.) take value from $1,2,\ldots$. Indices $j_1$, $s_1$ and their primes take value from $\{0,\ldots,L_1-1\}$. Indices $j_2$ and $s_2$ and their primes take value from $\{0,1,\ldots,L_2-1\}$. Also addition and subtraction for indices $j_1$, $s_1$ and their primes are performed in module $L_1$. Similarly, addition and subtraction for indices $j_2$, $s_2$ and their primes are performed in module $L_2$. 
\end{remark}

Note that 
\begin{align} f(x_1-s_1,x_2-s_2)&=
\sum_{i_1,i_2=1}^{\infty}\sum_{j_1=0,j_2=0}^{L_1-1,L_2-1} a_{j_1,j_2}^{i_1,i_2}\xi_{j_1}^{i_1}(x_1-s_1)\xi^{i_2}_{j_2}(x_2-s_2) \notag \\
&=\sum_{i_1,i_2=1}^{\infty}\sum_{j_1=0,j_2=0}^{L_1-1,L_2-1} a_{j_1,j_2}^{i_1,i_2}\xi_{j_1+s_1}^{i_1}(x_1)\xi^{i_2}_{j_2+s_2}(x_2) \notag \\
&= \sum_{i_1,i_2=1}^{\infty}\sum_{j_1=0,j_2=0}^{L_1-1,L_2-1} a_{j_1-s_1,j_2-s_2}^{i_1,i_2}\xi_{j_1}^{i_1}(x_1)\xi^{i_2}_{j_2}(x_2). \label{equation:shift2D}
\end{align}
Also, suppose 
\[  g(x_1,x_2)=\sum_{i_1,i_2=1}^{\infty}\sum_{j_1=0,j_2=0}^{L_1-1,L_2-1} b_{j_1,j_2}^{i_1,i_2}\xi_{j_1}^{i_1}(x_1)\xi_{j_2}^{i_2}(x_2). \]
Because SOBFs in 2D form an orthonormal set:
\begin{equation}
\int g(x_1,x_2)^*f(x_1-s_1,x_2-s_2) \dx= \sum_{i_1,i_2=1}^{\infty}\sum_{j_1=0,j_2=0}^{L_1-1,L_2-1}  {b_{j_1,j_2}^{i_1,i_2}}^* a_{j_1-s_1,j_2-s_2}^{i_1,i_2}.
\label{equation:shiftMultiplication}\end{equation}

We use a specific multi-index notation
\[ (i_1,i_2; j_1,j_2)_{i_1=1,i_2=1,j_1=0,j_2=0}^{i_1=\infty,i_2=\infty,j_1=L_1-1,j_2=L_2-1}, \] 
to refer to the entries of an infinite dimensional vector. To be rigorous, there is a (non-unique) one-to-one and onto mapping  
\[ \rho:\mathbb{N}\times \mathbb{N}\times \{0,\ldots,L_1-1\}\times\{0,\ldots,L_2-1\}\rightarrow \mathbb{N}. \] 
 Indeed, by $(i_1,i_2;j_1,j_2)$, we mean $\rho((i_1,i_2;j_1,j_2))$; however, to avoid cumbersome notation, we just use $(i_1,i_2;j_1,j_2)$ to refer to the positive integer instead.
As it will be seen shortly, $i_1$ and $j_1$ are related to depth index and shift index for the first coordinate, respectively. Similarly, $i_2$ and $j_2$ are related to the depth index and shift index for the second coordinate, respectively. We also use multi-index notation 
\[ (s_1,s_2)_{s_1=0,s_2=0}^{s_1=L_1-1,s_2=L_2-1},\]
to refer to the entries of a vector of length $L_1L_2$. Again, to be rigorous, there is a (non-unique) one-to-one and onto mapping  
\[ \tilde{\rho}:\{0,\ldots,L_1-1\}\times\{0,\ldots,L_2-1\}\rightarrow \{1,2,\ldots,L_1L_2\}. \] 
 Indeed, by $(s_1,s_2)$, we mean $\tilde{\rho}((s_1,s_2))$; however, to avoid cumbersome notation, we just use $(s_1,s_2)$ to refer to the elements of $\{1,\ldots,L_1L_2\}$.

In view of \eqref{equation:fSOBFexpansion2D}, function $f$ can be represented in SOBFs basis using infinite dimensional vector 
\begin{equation} \SOBF(f)(i_1,i_2;j_1,j_2):= a_{j_1,j_2}^{i_1,i_2}. \label{equation:definitionSOBFoperator}\end{equation}
Also for any infinite dimensional vector $a\in \mathbb{C}^{\mathbb{N}}$, define
\begin{equation} \ISOBF(a):= \sum_{i_1,i_2=1}^{\infty}\sum_{j_1=0,j_2=0}^{L_1-1,L_2-1} a(i_1,i_2;j_1,j_2)\xi_{j_1}^{i_1}(x_1)\xi_{j_2}^{i_2}(x_2).\label{equation:definitionISOBFoperator} \end{equation}

For any infinite dimensional vector $v$ and all pairs $(s_1,s_2)$, define transformation $S(s_1,s_2)$ in the following way:
\[ \tilde{v}=S(s_1,s_2)v  \qquad \hbox{if and only if} \qquad \tilde{v}(i_1,i_2;j_1,j_2)=v(i_1,i_2;j_1-s_1,j_2-s_2).\]
Note that from \eqref{equation:shift2D},
\[ \SOBF(f(x_1-s_1,x_2-s_2))=S(s_1,s_2)\SOBF(f(x_1,x_2)).\]
In view of \eqref{equation:shiftMultiplication}, for $s_1,s'_1=0,\ldots,L_1-1$ and $s_2,s'_2=0,\ldots,L_2-1$:
\begin{equation} \langle g(x_1-s_1,x_2-s_2),f(x_1-s'_1,x_2-s'_2) \rangle = \langle S(s_1,s_2)\SOBF(g),S(s'_1,s'_2)\SOBF(f) \rangle.
\label{equation:gfSS}\end{equation}
Define Set of Shift Orthogonal as follow: 
\begin{align} \mathcal{SSO}&=\{ v\in \mathbb{C}^{\mathbb{N}}: \langle v,S(s_1,s_2)v \rangle=\delta_{0s_1}\delta_{0s_2} \hbox{ for all $s_1$ and $s_2$} \}  \notag\\
&=\{ v\in \mathbb{C}^{\mathbb{N}}: \langle S(s'_1,s'_2) v,S(s_1,s_2)v \rangle=\delta_{s'_1s_1}\delta_{s'_2s_2} \hbox{ for all $s_1$, $s'_1$, $s_2$ and $s'_2$}\}.  \notag
 \end{align}
Observe that $\mathcal{SSO}$ is only a set and not a subspace. Moreover, equation \eqref{equation:gfSS} implies that $f(x_1,x_2)$ is shift orthogonal if and only if $\SOBF(f)\in \mathcal{SSO}$ (i.e. see theorem \ref{theorem:3equivalence2D}).

Define $\B$-transform to be operator $\B:\mathbb{C}^{\mathbb{N}}\rightarrow \mathbb{C}^{\mathbb{N}}$ defined by 
\begin{equation} \B(v)(i_1,i_2;j_1,j_2)=\sum_{\ell_1,\ell_2} e^{i2\pi(\frac{j_1}{L_1},\frac{j_2}{L_2})\cdot (\ell_1,\ell_2)}v(i_1,i_2;\ell_1,\ell_2). \label{equation:Bdefinition}\end{equation}
The intuition for $\B$-transform is that for every fixed $i_1$ and $i_2$, if $v(i_1,i_2;~:~,~:~)$ and $\B(v)(i_1,i_2;~:~,~:~)$ are thought as $L_1\times L_2$ matrices then 
\[ \B(v)(i_1,i_2;~:~,~:~)=L_1L_2\mathcal{F}_{2D}^{-1}(v(i_1,i_2;~:~,~:~), \]
where $\mathcal{F}_{2D}^{-1}$ is the 2D discrete inverse Fourier transform. The inverse of $\B$-transform, $\B^{-1}:\mathbb{C}^{\mathbb{N}}\rightarrow \mathbb{C}^{\mathbb{N}}$, is 
\[ \B^{-1}(v)(i_1,i_2;j_1,j_2)=\frac{1}{L_1L_2}\sum_{\ell_1,\ell_2} e^{-i2\pi(\frac{j_1}{L_1},\frac{j_2}{L_2})\cdot (\ell_1,\ell_2)}v(i_1,i_2;\ell_1,\ell_2). \]  
Similar to the above, $\B^{-1}$-transform can be written as 
\[ \B^{-1}(v)(i_1,i_2;~:~,~:~)=\frac{1}{L_1L_2}\mathcal{F}_{2D}(v(i_1,i_2;~:~,~:~)), \]
where $\mathcal{F}_{2D}$ is the 2D discrete Fourier transform.

The importance of $\B$-transform appears in the following two theorems:
\begin{theorem}
For given function $f$, the followings are equivalent:
\begin{enumerate}
\item $f$ is shift orthogonal. 
\item $\SOBF(f)\in \mathcal{SSO}$.
\item for all $(j_1,j_2)\in \{0,\ldots,L_1-1\}\times \{0,\ldots,L_2-1\}$,
\[ \|\B(\SOBF(f))(~:~,~:~;j_1,j_2) \|_2=1. \]
\end{enumerate}
\label{theorem:3equivalence2D}\end{theorem} 

\begin{theorem}
For given functions $f$ and $g$, the followings are equivalent:
\begin{enumerate}
\item  for all  $s_1,s'_1=0,\ldots,L_1-1$ and $s_2,s'_2=0,\ldots,L_2-1$,
\[ \langle g(x_1-s_1,x_2-s_2),f(x_1-s'_1,x_2-s'_2) \rangle =0. \]
\item  for all  $s_1=0,\ldots,L_1-1$ and $s_2=0,\ldots,L_2-1$,
\[ \langle \SOBF(g),S(s_1,s_2)\SOBF(f) \rangle =0. \]
\item for all  $j_1=0,\ldots,L_1-1$ and $j_2=0,\ldots,L_2-1$,
\[  \Big \langle \B(\SOBF(g))(~:~,~:~;j_1,j_2), \B(\SOBF(f))(~:~,~:~;j_1,j_2)\Big\rangle =0.\]
\end{enumerate}
\label{theorem:shiftPerpendicular2D}\end{theorem}
In order to prove the above two theorems, we need to introduce some more notations and concepts. We use multi-index notation to define matrices in the following way: Identify entries of matrix $A$ by $A(\cdot|\cdot)$, where the index to the left of ``$|$'' determines the row number of the entry and the index to the right of ``$|$'' determines the column number of the entry. Let $W$ be the $L_1L_2\times L_1L_2$ matrix defined by 
\[ W(s_1,s_2|j_1,j_2)=\frac{1}{\sqrt{L_1L_2}}e^{-i2\pi(\frac{s_1}{L_1},\frac{s_2}{L_2})\cdot(j_1,j_2)}. \]
Let $W_\infty$ to be infinite dimensional square matrix defined by 
\[ W_\infty(i_1,i_2;s_1,s_2|i'_1,i'_2;j_1,j_2)=\frac{1}{\sqrt{L_1L_2}}e^{-i2\pi(\frac{s_1}{L_1},\frac{s_2}{L_2})\cdot(j_1,j_2)}\delta_{i_1i'_1}\delta_{i_2i'_2}.\]
Finally, for any infinite dimensional vector $v$, let $\circulant(v)$ be the $L_1L_2\times \mathbb{N}$ matrix defined by
\[ \circulant(v)(s_1,s_2|i_1,i_2;j_1,j_2)=(S(s_1,s_2)v)(i_1,i_2;j_1,j_2)=v(i_1,i_2;j_1-s_1,j_2-s_2). \]
Note that $W$ and $W_\infty$ are unitary matrix (i.e. although the latter is infinite dimensional and by unitary we mean that $W_\infty W^\dagger_\infty$ is infinite dimensional diagonal matrix whose diagonal entries are one). To see this, observe that for example, 
\begin{align}
~& (W_\infty W^\dagger_\infty)(i_1,i_2;s_1,s_2|i'_1,i'_2;s'_1,s'_2) \notag\\
=& \sum_{i''_1,i''_2}\sum_{s''_1,s''_2} W_\infty(i_1,i_2;s_1,s_2|i''_1,i''_2;s''_1,s''_2)W^\dagger_\infty(i''_1,i''_2;s''_1,s''_2|i'_1,i'_2;s'_1,s'_2) \notag \\
=& \sum_{i''_1,i''_2}\sum_{s''_1,s''_2} W_\infty(i_1,i_2;s_1,s_2|i''_1,i''_2;s''_1,s''_2)\overline{W_\infty(i'_1,i'_2;s'_1,s'_2|i''_1,i''_2;s''_1,s''_2)} \notag \\
=&\sum_{i''_1,i''_2}\sum_{s''_1,s''_2} \frac{1}{\sqrt{L_1L_2}}e^{-i2\pi(\frac{s_1}{L_1},\frac{s_2}{L_2})\cdot(s''_1,s''_2)}\delta_{i_1i''_1}\delta_{i_2i''_2}\frac{1}{\sqrt{L_1L_2}}e^{i2\pi(\frac{s'_1}{L_1},\frac{s'_2}{L_2})\cdot(s''_1,s''_2)}\delta_{i'_1i''_1}\delta_{i'_2i''_2} \notag \\
=&\sum_{i''_1}\delta_{i_1i''_1}\delta_{i'_1i''_1}\sum_{i''_2}\delta_{i_2i''_2}\delta_{i'_2i''_2}\sum_{s''_1,s''_2}\frac{1}{L_1L_2}e^{i2\pi(\frac{s'_1-s_1}{L_1},\frac{s'_2-s_2}{L_2})\cdot(s''_1,s''_2)} \notag \\
=&\delta_{i_1i'_1}\delta_{i_2i'_2}\delta_{s_1s'_1}\delta_{s_2s'_2}. \notag
\end{align}

Another important property that $W$ and $W_\infty$ have is that if $\circulant(v)$ is multiplied on left by $W$ and on right by $W^{\dagger}_{\infty}$, then 
\begin{equation} (W\circulant(v)W^\dagger_\infty)(s_1,s_2|i_1,i_2;j_1,j_2)=\delta_{s_1j_1}\delta_{s_2j_2}\B(v)(i_1,i_2;j_1,j_2).
\label{equation:WcircW}\end{equation}
To see the above, note that 
\begin{align}
~& (W\circulant(v)W^\dagger_\infty)(s_1,s_2|i_1,i_2;j_1,j_2) \notag \\
= & \sum_{j'_1,j'_2}\sum_{i'_1,i'_2,s'_1,s'_2}W(s_1,s_2|j'_1,j'_2)\circulant(v)(j'_1,j'_2|i'_1,i'_2;s'_1,s'_2)W^{\dagger}_\infty(i'_1,i'_2;s'_1,s'_2|i_1,i_2;j_1,j_2)  \notag \\
= & \sum_{j'_1,j'_2}\sum_{i'_1,i'_2}\sum_{s'_1,s'_2}\frac{e^{-i2\pi(\frac{s_1}{L_1},\frac{s_2}{L_2})\cdot(j'_1,j'_2)}}{\sqrt{L_1L_2}}v(i'_1,i'_2;s'_1-j'_1,s'_2-j'_2)\frac{e^{i2\pi(\frac{j_1}{L_1},\frac{j_2}{L_2})\cdot(s'_1,s'_2)}}{\sqrt{L_1L_2}}\delta_{i_1i'_1}\delta_{i_2i'_2} \notag \\
=& \sum_{j'_1,j'_2}\sum_{\ell_1,\ell_2}\frac{e^{-i2\pi(\frac{s_1}{L_1},\frac{s_2}{L_2})\cdot(j'_1,j'_2)}}{\sqrt{L_1L_2}}v(i_1,i_2;\ell_1,\ell_2)\frac{e^{i2\pi(\frac{j_1}{L_1},\frac{j_2}{L_2})\cdot(\ell_1+j'_1,\ell_2+j'_2)}}{\sqrt{L_1L_2}} \notag \\
=& \sum_{j'_1,j'_2}\frac{1}{L_1L_2}e^{i2\pi(\frac{j_1-s_1}{L_1},\frac{j_2-s_2}{L_2})\cdot (j'_1,j'_2)}\sum_{\ell_1,\ell_2}v(i_1,i_2;\ell_1,\ell_2)e^{i2\pi(\frac{j_1}{L_1},\frac{j_2}{L_2})\cdot (\ell_1,\ell_2)} \notag \\
=& \delta_{s_1j_1}\delta_{s_2j_2}\B(v)(i_1,i_2;j_1,j_2), \notag
\end{align}
where \eqref{equation:Bdefinition} was used for the last equality. Equality \eqref{equation:WcircW} is very significant: it shows how $\circulant(v)$ can be turned into a ``pseudo-diagonal'' matrix using unitary matrices $W$ and $W_\infty$. Equation \eqref{equation:WcircW} is used extensively, in the remainder of this section.  

Finally, for any two infinite dimensional vectors $a$ and $b$:
\begin{equation} (\circulant(b)\circulant(a)^\dagger)(s_1,s_2|j_1,j_2)=\overline{\langle S(s_1,s_2)b,S(j_1,j_2)a \rangle}.
\label{equation:circacircb}\end{equation}
For observe that 
\begin{align}
~& (\circulant(b)\circulant(a)^\dagger)(s_1,s_2|j_1,j_2) \notag \\
=& \sum_{i_1,i_2, j'_1,j'_2}\circulant(b)(s_1,s_2|i_1,i_2; j'_1,j'_2)\circulant(a)^\dagger(i_1,i_2;j'_1,j'_2|j_1,j_2)	\notag \\
=&  \sum_{i_1,i_2, j'_1,j'_2}\circulant(b)(s_1,s_2|i_1,i_2; j'_1,j'_2)\overline{\circulant(a)(j_1,j_2|i_1,i_2;j'_1,j'_2)}  	\notag \\
=& \sum_{i_1,i_2, j'_1,j'_2}(S(s_1,s_2)b)(i_1,i_2;j'_1,j'_2)\overline{(S(j_1,j_2)a)(i_1,i_2;j'_1,j'_2)} \notag \\
=& \overline{\langle S(s_1,s_2)b,S(j_1,j_2)a \rangle}. \notag
\end{align}
Now we are ready to prove theorem \ref{theorem:3equivalence2D} and \ref{theorem:shiftPerpendicular2D}.

\textbf{Proof of theorem \ref{theorem:3equivalence2D}:}
The equivalence of 1 and 2 follows easily from \eqref{equation:gfSS} (i.e. with $g=f$) and the definition of $\mathcal{SSO}$. It remains to prove the equivalence between 2 and 3:

Set $v=\SOBF(f)$. Equation \eqref{equation:circacircb} (i.e. with $a=b=v$) implies that $v\in \mathcal{SSO}$ if and only if $ (\circulant(v)\circulant(v)^\dagger)(s_1,s_2|j_1,j_2)=\delta_{s_1j_1}\delta_{s_2j_2}$ for all $s_1$, $j_1$, $s_2$ and $j_2$ (i.e. matrix $\circulant(v)\circulant(v)^\dagger$ is the $L_1L_2\times L_1L_2$ identity matrix). 

Next let 
\[ V=W\circulant(v)W^\dagger_{\infty}. \]
Compute $VV^\dagger$ in two ways: On one hand, because $W_\infty$ is unitary 
\begin{equation} VV^\dagger=W\circulant(v)W^\dagger_{\infty}W_\infty \circulant(v)^\dagger W^\dagger=W\circulant(v)\circulant(v)^\dagger W^\dagger.
\label{equation:VVdagger1} \end{equation}
On the other hand, \eqref{equation:WcircW} yields that 
\begin{align}
~&VV^\dagger(s_1,s_2|j_1,j_2) \notag \\
=& \sum_{i_1,i_2, j'_1,j'_2} V(s_1,s_2|i_1,i_2;j'_1,j'_2)V^\dagger(i_1,i_2;j'_1,j'_2|j_1,j_2) \notag \\
=& \sum_{i_1,i_2, j'_1,j'_2}  V(s_1,s_2|i_1,i_2;j'_1,j'_2) \overline{V(j_1,j_2|i_1,i_2;j'_1,j'_2)} \notag \\
=& \sum_{i_1,i_2, j'_1,j'_2}  \B(v)(i_1,i_2;j'_1,j'_2)\delta_{s_1j'_1}\delta_{s_2j'_2}\overline{\B(v)(i_1,i_2;j'_1,j'_2)\delta_{j_1j'_1}\delta_{j_2j'_2}} \notag \\
=& \sum_{i_1,i_2} |\B(v)(i_1,i_2; j_1,j_2)|^2\delta_{s_1j_1}\delta_{s_2j_2} \notag \\
=& \|\B(v)(~:~,~:~;j_1,j_2)\|_2^2 ~ \delta_{s_1j_1}\delta_{s_2j_2}. \label{equation:VVdagger2} 
\end{align}
Now if $ \circulant(v)\circulant(v)^\dagger$ is the identity matrix (i.e. $v \in \mathcal{SSO}$), then from \eqref{equation:VVdagger1} and $W$ being unitary, one concludes that $VV^\dagger$ is the identity matrix. Hence, by \eqref{equation:VVdagger2}, it must be the case that 
\[ \|\B(v)(~:~,~:~;j_1,j_2)\|_2=1, \]
for all $(j_1,j_2)\in \{0,\ldots,L_1-1\}\times\{0,\ldots,L_2-1\}$.

Conversely, if the above holds, then by \eqref{equation:VVdagger2}, $VV^\dagger$ is the identity matrix. Therefore, because $W$ is a unitary matrix, $W^\dagger VV^\dagger W$ would be the identity matrix as well. Equation \eqref{equation:VVdagger1}, yields that 
\[ W^\dagger VV^\dagger W=\circulant(v)\circulant(v)^\dagger. \] 
Hence, $\circulant(v)\circulant(v)^\dagger$ is the identity matrix, which implies $v\in \mathcal{SSO}$.
\qed 

\textbf{Proof of theorem \ref{theorem:shiftPerpendicular2D}:}
The equivalence of 1 and 2 follows easily from \eqref{equation:gfSS} and the fact that 
\[ \langle \SOBF(g),S(s_1,s_2)\SOBF(f) \rangle =0 \quad \hbox{ for all $s_1$ and $s_2$}, \]
is equivalent to 
\begin{equation} \langle S(s_1,s_2)\SOBF(g),S(s'_1,s'_2)\SOBF(f) \rangle =0 \quad \hbox{ for all $s_1$, $s'_1$, $s_2$ and $s'_2$}. \label{equation:alternative2}\end{equation}
It remains to show that \eqref{equation:alternative2} is equivalent to 3. Let $b=\SOBF(g)$ and $a=\SOBF(f)$. Equation \eqref{equation:circacircb}, yields that \eqref{equation:alternative2} is equivalent to $(\circulant(b)\circulant(a)^\dagger)(s_1,s_2|j_1,j_2)=0$ for all  $s_1$, $j_1$, $s_2$ and $j_2$ (i.e. matrix $\circulant(b)\circulant(a)^\dagger$ is the $L_1L_2\times L_1L_2$ zero matrix).

Next let 
\[ B=W\circulant(b)W^\dagger_{\infty} \qquad \hbox{and} \qquad A=W\circulant(a)W^\dagger_{\infty} . \]
Compute $BA^\dagger$ in two ways: On one hand, because $W_\infty$ is unitary 
\begin{equation} BA^\dagger=W\circulant(b)W^\dagger_{\infty}W_\infty \circulant(a)^\dagger W^\dagger=W\circulant(b)\circulant(a)^\dagger W^\dagger.
\label{equation:BAdagger1} \end{equation}
On the other hand, \eqref{equation:WcircW} yields that 
\begin{align}
~&BA^\dagger(s_1,s_2|j_1,j_2) \notag \\
=& \sum_{i_1,i_2, j'_1,j'_2} B(s_1,s_2|i_1,i_2;j'_1,j'_2)A^\dagger(i_1,i_2;j'_1,j'_2|j_1,j_2) \notag \\
=& \sum_{i_1,i_2, j'_1,j'_2}  B(s_1,s_2|i_1,i_2;j'_1,j'_2) \overline{A(j_1,j_2|i_1,i_2;j'_1,j'_2)} \notag \\
=& \sum_{i_1,i_2, j'_1,j'_2}  \B(b)(i_1,i_2;j'_1,j'_2)\delta_{s_1j'_1}\delta_{s_2j'_2}\overline{\B(a)(i_1,i_2;j'_1,j'_2)\delta_{j_1j'_1}\delta_{j_2j'_2}} \notag \\
=& \sum_{i_1,i_2} \B(b)(i_1,i_2; j_1,j_2)\overline{\B(a)(i_1,i_2; j_1,j_2)}\delta_{s_1j_1}\delta_{s_2j_2} \notag \\
=& \overline{\langle\B(b)(~:~,~:~;j_1,j_2),\B(a)(~:~,~:~;j_1,j_2)\rangle} ~ \delta_{s_1j_1}\delta_{s_2j_2}. \label{equation:BAdagger2} 
\end{align}
Now if $ \circulant(b)\circulant(a)^\dagger$ is the zero matrix, then \eqref{equation:BAdagger1} implies that $BA^\dagger$ is the zero matrix. Hence, by \eqref{equation:BAdagger2}, it must be the case that 
\[ \langle \B(b)(~:~,~:~;j_1,j_2),\B(a)(~:~,~:~;j_1,j_2) \rangle=0,  \]
for all $(j_1,j_2)\in \{0,\ldots,L_1-1\}\times\{0,\ldots,L_2-1\}$. 

Conversely, if the above holds, then by \eqref{equation:BAdagger2}, $BA^\dagger$ is the zero matrix. Therefore, $W^\dagger BA^\dagger W$ would also be the zero matrix. Equation \eqref{equation:BAdagger1} yields that 
\[ W^\dagger BA^\dagger W=\circulant(b)\circulant(a)^\dagger. \] 
Hence, $\circulant(b)\circulant(a)^\dagger$ is the zero matrix, which implies \eqref{equation:alternative2}.
\qed

For any function $g\in L^2(\Omega)$, let $\Pi g$ denote the projection of $g$ into the set of shift orthogonal functions; that is, 
\[ \Pi g:= \argmin_f \|g-f\|_2 \qquad \hbox{ subject to $f$ being shift orthogonal. } \]
Indeed using SOBFs basis,
\[ \Pi g=\ISOBF(Proj_{\mathcal{SSO}}(\SOBF(g))), \]
where for any $b\in \mathbb{C}^{\mathbb{N}}$,
\[ Proj_{\mathcal{SSO}} (b) := \argmin_v \| b-v\|_2 \qquad \hbox{subject to } \qquad v\in \mathcal{SSO}.\]
Observe that in the above two definitions the minimum arguments are not necessarily unique, and $\Pi g$ and $Proj_{\mathcal{SSO}}(b)$ are sets.

Define operator $\Theta:\mathbb{C}^{\mathbb{N}}\rightarrow \mathbb{C}^{\mathbb{N}}$ by 
\[ \Theta(v)(:,:,j_1,j_2)=\begin{cases}
\vec{e} \qquad &\hbox{if } v(:,:,j_1,j_2)=\vec{0} \\
\frac{v(:,:,j_1,j_2)}{\|v(~:~,~:~,j_1,j_2)\|_2} \qquad &\hbox{otherwise},
\end{cases} \] 
where $\vec{e}$ is a fixed infinite dimensional real vector with unit $L^2$ norm.
\begin{lemma}
For any $b\in \mathbb{C}^{\mathbb{N}}$,
\[   \B^{-1}(\Theta(\B(b))) \in Proj_{\mathcal{SSO}}(b).\]
\label{lemma:projection2D}\end{lemma}
\proof Suppose $v\in \mathcal{SSO}$ and set 
%\begin{equation} B=W\circulant(b)W^\dagger_\infty \qquad \hbox{ and } \qquad V=W\circulant(v)W^\dagger_\infty. \label{equation:BA} \end{equation}
\[ p=\B(b) \qquad \hbox{ and } q=\B(v). \]

Note that,
\begin{align}
 \|b-v\|_2^2=\frac{1}{L_1L_2}\|\circulant(b)-\circulant(v)\|_F^2& =\frac{1}{L_1L_2}\|W\circulant(b) W^\dagger_\infty-W\circulant(v)W^\dagger_\infty\|_F^2   \notag \\
&=\frac{1}{L_1L_2}\|\B(b)-\B(v)\|_2^2=\frac{1}{L_1L_2}\|p-q\|_2^2, \notag
\end{align}
where equalities similar to \eqref{equation:WcircW} were used for the second last inequality. Now minimizing $\|p-q\|_2$ amounts to solving $L_1L_2$ subproblems: for every $j_1=0,\ldots,L_1-1$ and $j_2=0,\ldots,L_2-1$, solve
\begin{equation} \argmin \sum_{i_1,i_2} |p(i_1,i_2;j_1,j_2)-q(i_1,i_2;j_1,j_2)|^2 \quad \hbox{subject to }  \sum_{i_1,i_2} |q(i_1,i_2;j_1,j_2)|^2=1 \label{equation:subProblem2D}. \end{equation}
The constraints in the above subproblems are due to $v\in \mathcal{SSO}$ and equivalence of 2 and 3 in theorem \ref{theorem:3equivalence2D}. 
The solutions to the above subproblems are exactly
\[  \begin{cases}
q(:,:;j_1,j_2)=p(:,:;j_1,j_2)/\|p(:,:;j_1,j_2)\|_2  \quad &\hbox{ if } \|p(:,:;j_1,j_2)\|_2\neq 0  \\
\hbox{any infinite dimensional complex vector with unit $L^2$ norm} \quad &\hbox{ if } \|p(:,:;j_1,j_2)\|_2= 0,
\end{cases}\]
that is, projection of $p(:,:;j_1,j_2)$ into an infinite dimensional ball of radius 1. In particular, set $q=\Theta(p)$ (i.e. choose a fixed real valued vector in the second case above), in which case, $v=\B^{-1}(\Theta(p))$ would be an element of $Proj_{\mathcal{SSO}}(b)$. \qed

\begin{remark}\label{remark:realness}
Theorem 2.2 in paper \cite{BIB:Chu} and definition of operator $\Theta$, yields that if infinite dimensional vector $b$ is real valued, then $\B^{-1}(\Theta(\B(b)))$ is also real valued. This is why in the case $ v(:,:,j_1,j_2)=\vec{0}$, we assign to operator $\Theta$ the fixed infinite dimensional vector $\vec{e}$; which is real valued and has unit $L^2$ norm. Indeed, (as it is apparent in the proof of lemma \ref{lemma:projection2D}) had we defined operator $\Theta$ to output all infinite dimensional (complex) vector of unit $L^2$ norm in the case $ v(:,:,j_1,j_2)=\vec{0}$, then $\B^{-1}(\Theta(\B(b)))$ would have been a set equal to $ Proj_{\mathcal{SSO}}(b)$; however, some elements of $\B^{-1}(\Theta(\B(b)))$ would have been complex valued.
\end{remark}

\section{Fast Algorithm for Projection to the Space of Shift Orthogonal Functions}\label{section:projectionAlgorithm}
For computational purposes, only a finite number of SOBFs basis are used to represent a function. In this section, assume that for the first coordinate all SOBFs basis whose depth index is smaller or equal to $N_1$ and for the second coordinate all SOBFs basis whose depth index is smaller or equal to $N_2$ are used to denote functions in $L^2(\Omega)$. That is,
\[ g(x_1,x_2)=\sum_{i_1=1,i_2=1}^{\infty}\sum_{j_1=0,j_2=0}^{L_1-1,L_2-2} b_{j_1,j_2}^{i_1,i_2}\xi_{j_1}^{i_1}(x_1)\xi_{j_2}^{i_2}(x_2) \approx\sum_{i_1,i_2=1}^{N_1,N_2}\sum_{j_1=0,j_2=0}^{L_1-1,L_2-1} b_{j_1,j_2}^{i_1,i_2}\xi_{j_1}^{i_1}(x_1)\xi_{j_2}^{i_2}(x_2). \]

The analysis done in section \ref{section:theorySOBF2D} can be adapted for this situation by simple modification. In particular, index $i_1$ (and $i'_1$) takes value from $1,\ldots,N_1$ instead of $1,2,\ldots$ and index $i_2$ (and $i'_2$) takes value from $1,\ldots,N_2$ instead of $1,2,\ldots$.   

The adapted definition for set $\mathcal{SSO}$ with finite depth indices is 
\[ \mathcal{SSO}(N_1N_2)=
\{ v\in \mathbb{C}^{N_1N_2L_1L_2}: \langle v,S(s_1,s_2)v \rangle=\delta_{0s_1}\delta_{0s_2} \hbox{ for all $s_1$ and $s_2$}\}. \]

As noted earlier, an important question that arises in optimization problems that involve shift orthogonality constraints is to find $\Pi g$ for a given function $g$; that is, find shift orthogonal function $f$ that minimizes $\|g-f\|_2$. When functions are expressed in terms of tensor product of one dimensional SOBFs basis (i.e. with corresponding depth indices smaller or equal to $N_1$ and $N_2$), then the question is equivalent to: given $\vec{b}\in \mathbb{C}^{N_1N_2L_1L_2}$, solve 
\begin{equation} 
Proj_{\mathcal{SSO}(N_1N_2)} ( \vec{b} ) = \argmin \|\vec{b}-\vec{v}\|_2 \qquad \hbox{subject to } \qquad \vec{v}\in \mathcal{SSO}(N_1N_2).
\label{equation:SORoptimization} 
\end{equation}

Result of lemma \ref{lemma:projection2D} in section \ref{section:theorySOBF2D} implies that the solution to problem \eqref{equation:SORoptimization} can be obtained using the procedure in algorithm \ref{algorithm:procedureN2D}.
 
\begin{algorithm}[H]
\SetKwInOut{Input}{input}\SetKwInOut{Output}{output}
\SetKwComment{Comment}{}{}
\KwIn{$\vec{b}$}
\KwOut{$\vec{v} = Proj_{\mathcal{SSO}(N_1N_2)} ( \vec{b} ) $}
\BlankLine
%\dontprintsemicolon
\For{$i_1=1,\ldots,N_1$ and $i_2=1,\ldots,N_2$}{
$p(i_1,i_2;~:~,~:~)=L_1L_2 \mathcal{F}_{2D}^{-1}(b(i_1,i_2;~:~,~:~))$   \tcp*[r]{$p=\B(b).$} 
}
%Set $B=W^\dagger \circulant_N(b) W_N$. \;
\For{$j_1=0,\ldots,L_1-1$ and $j_2=0,\ldots,L_2-1$}{
%Set $P_i=\sqrt{\sum_{k=1}^{N} |B_{i,(k-1)L+i}|^2}$\;
\eIf {$\|p(:,:;j_1,j_2)\|_2\neq 0$}{
	$q(:,:;j_1,j_2)=p(:,:;j_1,j_2)/\|p(:,:;j_1,j_2)\|_2$  \tcp*[r]{$q=\B(v).$}
}{
	$q(:,:;j_1,j_2)=\vec{\mathbf 1}/\sqrt{N_1N_2}$.
}
}
\For{$i_1=1,\ldots,N_1$ and $i_2=1,\ldots,N_2$}{
$ v(i_1,i_2;~:~,~:~)=\frac{1}{L_1L_2} \mathcal{F}_{2D}(q(i_1,i_2;~:~,~:~))$  \tcp*[r]{$v=\B^{-1}(q).$}
}
%$a$ is the first row of matrix $WAW_N^\dagger$.
\caption{Projection to $\mathcal{SSO}(N_1N_2)$}
\label{algorithm:procedureN2D}
\end{algorithm}

All the results that were developed in section \ref{section:theorySOBF2D} can also be easily adapted for domains with other dimensions. For example suppose $\Omega=[0,L_1]\times [0,L_2]\times [0,L_3]$ (i.e. using appropriate scaling, it is assumed that the length of the shift along each coordinate is 1), and let 
\[g(x_1,x_2,x_3)\approx \sum_{i_1=1,i_2=1,i_3=1}^{N_1,N_2,N_3}\sum_{j_1=0,j_2=0,j_3=0}^{L_1-1,L_2-1,L_3-1} b_{j_1,j_2,j_3}^{i_1,i_2,i_3}\xi_{j_1}^{i_1}(x_1)\xi_{j_2}^{i_2}(x_2)\xi_{j_3}^{i_3}(x_3).  \]
Three dimensional version of algorithm \ref{algorithm:procedureN2D} is:

\begin{algorithm}[H]
\SetKwInOut{Input}{input}\SetKwInOut{Output}{output}
\SetKwComment{Comment}{}{}
\KwIn{$\vec{b}$}
\KwOut{$\vec{v} = Proj_{\mathcal{SSO}(N_1N_2N_3)} ( \vec{b} ) $}
\BlankLine
%\dontprintsemicolon
\For{$i_1=1,\ldots,N_1$, $i_2=1,\ldots,N_2$ and  $i_3=1,\ldots,N_3$}{
$p(i_1,i_2,i_3;~:~,~:~,~:~)=L_1L_2L_3 \mathcal{F}_{3D}^{-1}(b(i_1,i_2,i_3;~:~,~:~,~:~))$   \tcp*[r]{$p=\B(b).$} 
}
%Set $B=W^\dagger \circulant_N(b) W_N$. \;
\For{$j_1=0,\ldots,L_1-1$, $j_2=0,\ldots,L_2-1$ and $j_3=0,\ldots,L_3-1$}{
%Set $P_i=\sqrt{\sum_{k=1}^{N} |B_{i,(k-1)L+i}|^2}$\;
\eIf {$\|p(:,:,:;j_1,j_2,j_3)\|_2\neq 0$}{
	$q(:,:,:;j_1,j_2,j_3)=p(:,:,:;j_1,j_2,j_3)/\|p(:,:,:;j_1,j_2,j_3)\|_2$  \tcp*[r]{$q=\B(v).$}
}{
	$q(:,:,:;j_1,j_2,j_3)=\vec{\mathbf 1}/\sqrt{N_1N_2N_3}$.
}
}
\For{$i_1=1,\ldots,N_1$, $i_2=1,\ldots,N_2$ and $i_3=1,\ldots,N_3$}{
$ v(i_1,i_2,i_3;~:~,~:~,~:~)=\frac{1}{L_1L_2L_3} \mathcal{F}_{3D}(q(i_1,i_2,i_3;~:~,~:~,~:~))$  \tcp*[r]{$v=\B^{-1}(q).$}
}
%$a$ is the first row of matrix $WAW_N^\dagger$.
\caption{Projection to $\mathcal{SSO}(N_1N_2N_3)$}
\label{algorithm:procedureN3D}
\end{algorithm}

Finally, the one dimensional version of algorithm \ref{algorithm:procedureN2D} for domain $\Omega=[0,L]$ (i.e. again using appropriate scaling, it is assumed that the length of the shift is 1) and 
\[ g(x)\approx \sum_{i=1}^{N}\sum_{j=0}^{L-1}b_j^i \xi_j^i(x), \]
is the following:

\begin{algorithm}[H]
\SetKwInOut{Input}{input}\SetKwInOut{Output}{output}
\SetKwComment{Comment}{}{}
\KwIn{$\vec{b}$}
\KwOut{$\vec{v} = Proj_{\mathcal{SSO}(N)} ( \vec{b} ) $}
\BlankLine
%\dontprintsemicolon
\For{$i=1,\ldots,N$}{
$p(i;~:~)=L\mathcal{F}_{1D}^{-1}(b(i;~:~))$   \tcp*[r]{$p=\B(b).$} 
}
%Set $B=W^\dagger \circulant_N(b) W_N$. \;
\For{$j=0,\ldots,L-1$}{
%Set $P_i=\sqrt{\sum_{k=1}^{N} |B_{i,(k-1)L+i}|^2}$\;
\eIf {$\|p(~:~;j)\|_2\neq 0$}{
	$q(~:~;j)=p(~:~;j)/\|p(~:~;j)\|_2$  \tcp*[r]{$q=\B(v).$}
}{
	$q(~:~;j)=\vec{\mathbf 1}/\sqrt{N}$.
}
}
\For{$i=1,\ldots,N$}{
$ v(i;~:~)=\frac{1}{L} \mathcal{F}_{1D}(q(i;~:~))$  \tcp*[r]{$v=\B^{-1}(q).$}
}
%$a$ is the first row of matrix $WAW_N^\dagger$.
\caption{Projection to $\mathcal{SSO}(N)$}
\label{algorithm:procedureN1D}
\end{algorithm}

%\begin{figure}[h]
%\centering
%\begin{minipage}{0.8\linewidth}
%\includegraphics[width=1\linewidth]{figures/absProjection.eps}
%\end{minipage}\hfill\\  
%\begin{minipage}{0.8\linewidth}
%\centering
%\includegraphics[width=1\linewidth]{figures/squareProjection.eps}
%\end{minipage}\hfill\\
%\begin{minipage}{0.8\linewidth}
%\centering
%\includegraphics[width=1\linewidth]{figures/sinProjection.eps}
%\end{minipage}\hfill\\ 
%\caption{Application of algorithm \ref{algorithm:procedureN1D} to find $\Pi g$, projection of function $g$ to the set of shift orthogonal functions. From top to bottom, $g=|x-4|/4$, $g=(x-4)^2/16$ and $g=\sin(\pi x/2)$. Here $L=8$, $w=2$ (the length of the shifts) and $N=6$; that is we are using total of 24 SOPWs to represent functions.}
%\label{figure:projectionExample}
%\end{figure}

\subsection{Computational Complexity and Important Features of the Algorithm}\label{section:complexity}
This sections describes the computational complexity of algorithm \ref{algorithm:procedureN2D} and highlights some of the important properties of this algorithm. The results stay the same for domains with dimensions other than $d=2$.

Let $M=L_1L_2N_1N_2$ be the size of input vector $\vec{b}$; which indicates the number of coefficients used to represent the given function. Algorithm \ref{algorithm:procedureN2D} consists of three ``for'' loops. Each iteration in the first and the last ``for" loop can be computed using $O(L_1L_2\log(L_1L_2))$ operations via inverse Fast Fourier Transform and Fast Fourier Transform, respectively.  Each iteration in the second ``loop" can be done using $O(N_1N_2)$. Therefore, algorithm \ref{algorithm:procedureN2D} can be performed using 
\[ N_1N_2O(L_1L_2\log(L_1L_2))+L_1L_2O(N_1N_2)+N_1N_2O(L_1L_2\log(L_1L_2)) \]
operations, which leads to computational complexity of 
\[ O(M\log(L_1L_2))<O(M\log(M)). \]

Furthermore, note that each of the ``for" loops in algorithms \ref{algorithm:procedureN2D} can be done in parallel. This enhances the speed of the algorithm even further and makes it suitable for inputs with large dimensions.

Another nice property of algorithm \ref{algorithm:procedureN2D} is that for real valued input vector $\vec{v}$, it outputs a real valued vector $Proj_{\mathcal{SSO}(N_1N_2)} ( \vec{b} )$ (i.e. recall remark \ref{remark:realness}). The importance of this property is that if function $g$ and SOBFs $\{\xi_{j_1}^{i_1}(x_1)\xi_{j_2}^{i_2}(x_2)\}$ are real valued then $\SOBF(g)$ would be a real valued vector, and therefore using algorithm \ref{algorithm:procedureN2D}, the projected $\Pi g$ would also be real valued.

\section{An interesting example of SOBFs: Shift Orthogonal Plane Waves (SOPWs)}\label{section:definitionSOPW}
This section provides an example of real valued SOBFs with certain nice properties called Shift Orthogonal Plane Waves (SOPWs). As it will be seen shortly, SOPWs are suitable for numerical computation because there exist an exact prescription of them in terms of Fourier basis. Therefore, functions can be expanded in terms of SOPWs very efficiently using FFT and its inverse. 

Consider 1D domain $\Omega=[0,L]$ with periodic boundary and by scaling (i.e. replacing $L$ by $L/w$) assume that $w=1$. Furthermore, suppose that $L$ is even. This assumption is made so that the formulas provided in this section are easier to express. Nevertheless, the assumption that $L$ is even is not very restrictive as the parity of $L$ is not significant in many applications. 

Recall that functions 
\begin{equation} \phi_n(x)=\frac{1}{\sqrt{L}}e^{i2\pi n x/L} \qquad \hbox{ for } n\in \mathbb{Z}, \label{equation:exponentialBasis}\end{equation}
as well as 
\begin{equation}  \{\frac{1}{\sqrt{L}},\frac{2}{\sqrt{2L}}\cos(2\pi nx/L), \frac{2}{\sqrt{2L}}\sin(2\pi nx/L)\}_{n=1}^{\infty}\label{equation:sinCosBasis} \end{equation}
form orthonormal basis for $L^2(\Omega)$. Denote the $L$-th roots of unity by 
\[ \omega_j=e^{i2\pi j/L} \qquad \hbox{for } \qquad j=0, \ldots,L-1. \]
Shift Orthogonal Plane Waves (SOPWs) are denoted by 
\[ \{\theta^i_j(x)\}_{i=1,j=0}^{i=\infty,j=L-1},  \]
and defined in the following way: set
 \begin{align}
\theta^1_j(x)&=\frac{1}{\sqrt{L}}\sum_{|n|<\frac{L}{2}}\omega_j^{-n}\phi_n(x)+\frac{1}{\sqrt{2L}}\sum_{|n|=\frac{L}{2}}\omega_j^{-n}\phi_{n}(x) \label{equation:SOPWfourierLine1}  \\
\theta^k_j(x)&=\frac{1}{\sqrt{L}}\sum_{\frac{(k-1)L}{2}<|n|<\frac{kL}{2}}(\sgn(n)i)^{k-1}\omega_j^{-n}\phi_n(x)+\frac{1}{\sqrt{2L}}\sum_{|n|=\frac{(k-1)L}{2},\frac{kL}{2}}(\sgn(n)i)^{k-1}\omega_j^{-n}\phi_n(x) \label{equation:SOPWfourierLine2} 
\end{align}

\begin{figure}[h]
\centering
\begin{minipage}{0.8\linewidth}
\includegraphics[width=1\linewidth]{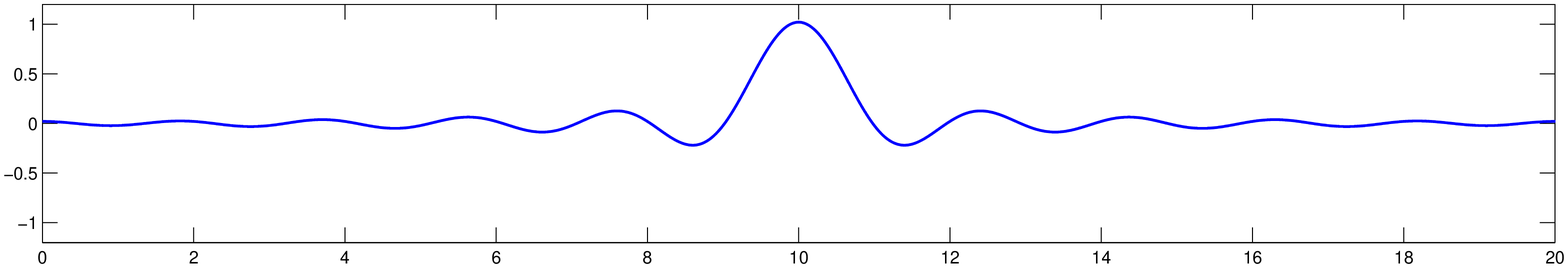}
\end{minipage}\hfill\\  
\begin{minipage}{0.8\linewidth}
\centering
\includegraphics[width=1\linewidth]{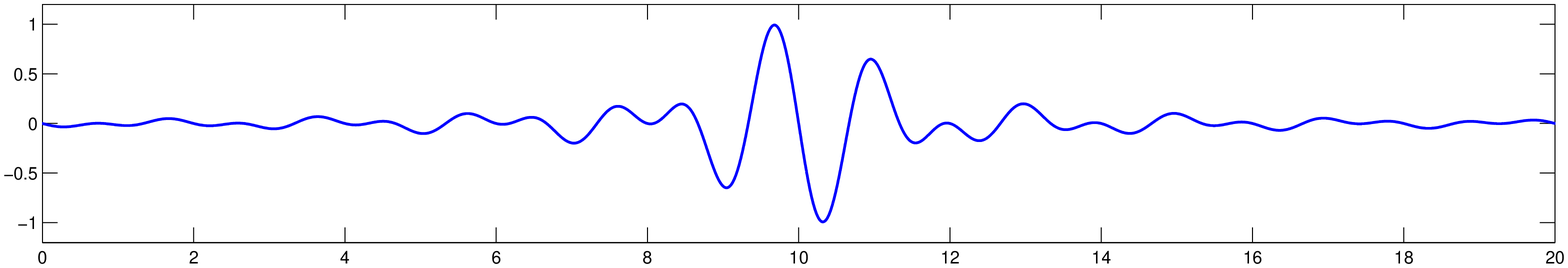}
\end{minipage}\hfill\\
\begin{minipage}{0.8\linewidth}
\centering
\includegraphics[width=1\linewidth]{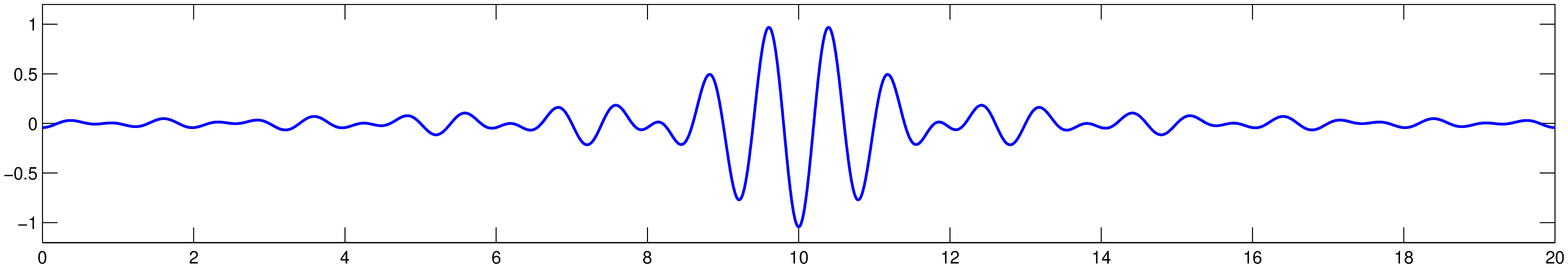}
\end{minipage}\hfill\\ 
\begin{minipage}{0.8\linewidth}
\centering
\includegraphics[width=1\linewidth]{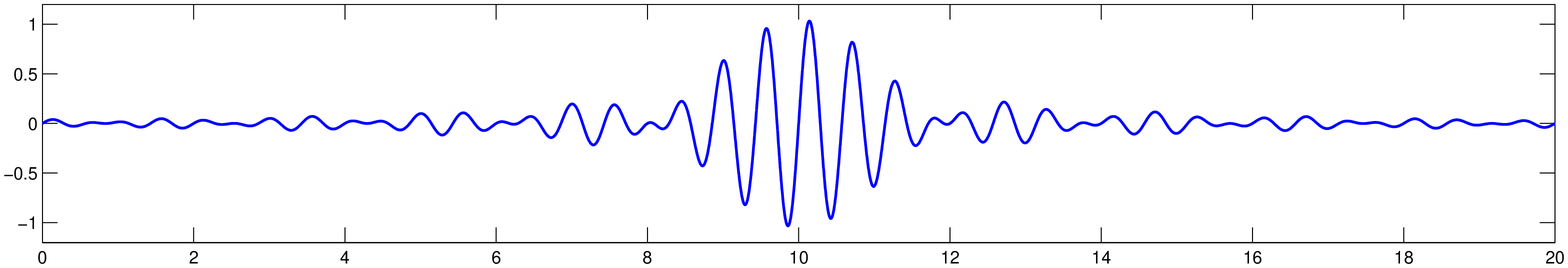}
\end{minipage}\hfill \\
\begin{minipage}{0.8\linewidth}
\centering
\includegraphics[width=1\linewidth]{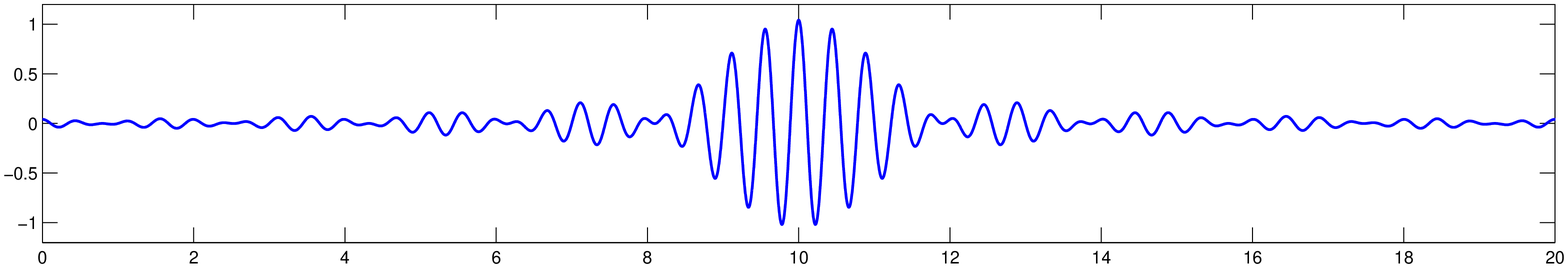}
\end{minipage}\hfill\\ 
\begin{minipage}{0.8\linewidth}
\centering
\includegraphics[width=1\linewidth]{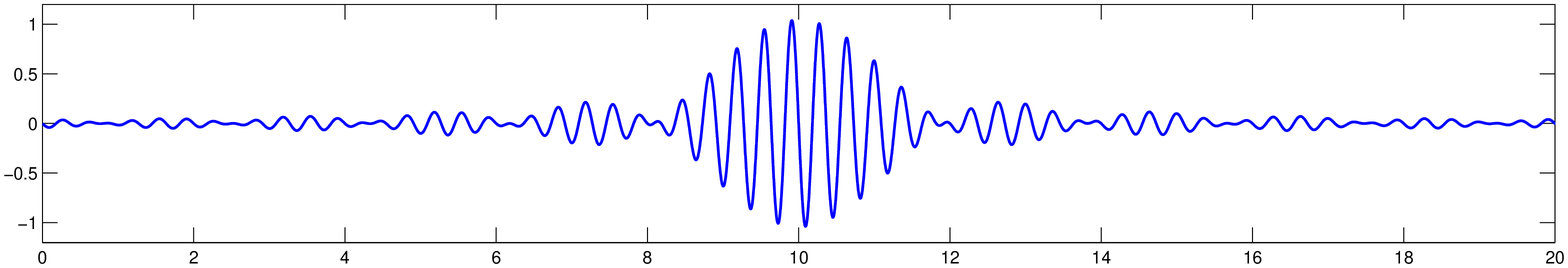}
\end{minipage}\hfill 
\caption{From top to bottom, the first 6 SOPWs given by equations \eqref{equation:SOPWfourierLine1} and \eqref{equation:SOPWfourierLine2} with distinct depth index and shift index equal to $L/2$.}
\label{figure:first_6}
\end{figure}

% \begin{align}
%\theta^1_j(x)&=\frac{1}{\sqrt{L}}\sum_{|n|<\frac{L}{2}}e^{-i2\pi nj/L}\phi_n(x)+\frac{1}{\sqrt{2L}}\sum_{|n|=\frac{L}{2}}e^{-i2\pi nj/L}\phi_{n}(x) \notag \\
%\theta^k_j(x)&=\frac{1}{\sqrt{L}}\sum_{\frac{(k-1)L}{2}<|n|<\frac{kL}{2}}(\sgn(n)i)^{k-1}e^{-i2\pi nj/L}\phi_n(x)+\frac{1}{\sqrt{2L}}\sum_{|n|=\frac{(k-1)L}{2},\frac{kL}{2}}(\sgn(n)i)^{k-1}e^{-i2\pi nj/L}\phi_n(x) \label{equation:SOPWfourier} 
%\end{align}
Figure \ref{figure:first_6} plots SOPWs given by equations \eqref{equation:SOPWfourierLine1} and \eqref{equation:SOPWfourierLine2} with depth index ranging from 1 to 6 and shift index equal to $L/2$.

Using the expression 
\begin{equation} \frac{1}{L} \sum_{n=k}^{k+L-1}\omega_j^n=\delta_{j0},  \label{equation:importantIdentity}\end{equation}
%\begin{equation} \frac{1}{L} \sum_{n=k}^{k+L-1}e^{i2\pi jn/L}=\delta_{j0},  \label{equation:importantIdentity}\end{equation}
and after some calculations one can verify that for $n\geq 1$,
\begin{align}
\phi_0(x)&=\frac{1}{\sqrt{L}}\sum_{j=0}^{L-1} \theta^1_j(x), \label{equation:fourierSOPWLine1}\\
\phi_n(x)&=\frac{(-\sgn(n) i)^{k-1}}{\sqrt{L}}\sum_{j=0}^{L-1}\omega_j^n \theta^k_j(x) \quad \hbox{for}\quad \frac{(k-1)L}{2}<|n|<\frac{kL}{2},\label{equation:fourierSOPWLine2}\\
\phi_n(x)&=\frac{(-\sgn(n) i)^{k-1}}{\sqrt{2L}}\left (\sum_{j=0}^{L-1}\omega_j^n\theta^k_j(x)- \sgn(n) i \sum_{j=0}^{L-1}\omega_j^n\theta^{k+1}_j(x) \right)  \hbox{ for } |n|=\frac{kL}{2}. \label{equation:fourierSOPWLine3}
\end{align} 
%{\small \begin{align}
%\phi_0(x)&=\frac{1}{\sqrt{L}}\sum_{j=0}^{L-1} \theta^1_j(x), \label{equation:fourierSOPWLine1}\\
%\phi_n(x)&=\frac{(-\sgn(n) i)^{k-1}}{\sqrt{L}}\sum_{j=0}^{L-1}e^{i2\pi nj/L} \theta^k_j(x) \quad \hbox{for}\quad \frac{(k-1)L}{2}<|n|<\frac{kL}{2},\label{equation:fourierSOPWLine2}\\
%\phi_n(x)&=\frac{(-\sgn(n) i)^{k-1}}{\sqrt{2L}}\left (\sum_{j=0}^{L-1}e^{i2\pi nj/L} \theta^k_j(x)- \sgn(n) i \sum_{j=0}^{L-1}e^{i2\pi nj/L} \theta^{k+1}_j(x) \right)  \hbox{ for } |n|=\frac{kL}{2}. \label{equation:fourierSOPWLine3}
%\end{align} }

First note that $\theta_j^i(x)=\theta_0^i(x-j)$. Moreover, it is straightforward using identity \eqref{equation:importantIdentity} to verify that $\{\theta^i_j\}_{i=1,j=0}^{i=\infty,j=L-1}$ form an orthonormal set. Finally, $\{\theta^i_j\}_{i=1,j=0}^{i=\infty,j=L-1}$ is complete in $L^2(\Omega)$ because of relations \eqref{equation:fourierSOPWLine1}, \eqref{equation:fourierSOPWLine2}, \eqref{equation:fourierSOPWLine3} and completeness of $\{\phi_n\}_{n=-\infty}^{n=\infty}$. Hence, the set of SOPWs defined by \eqref{equation:SOPWfourierLine1} and \eqref{equation:SOPWfourierLine2} is an example of SOBFs.

Equations \eqref{equation:SOPWfourierLine1} and \eqref{equation:SOPWfourierLine2} can be re-written to represent SOPWs $\{\theta^i_j\}_{i=1,j=0}^{i=\infty,j=L-1}$ in terms of \eqref{equation:sinCosBasis}:
{\footnotesize \begin{align}
\theta^1_j(x) &= \frac{1}{L}+\sum_{n=1}^{L/2-1}\frac{2}{L}\cos(\frac{2\pi n(x-j)}{L})+\frac{\sqrt{2}}{L}\cos(\frac{2\pi(L/2)(x-j)}{L}),  \notag\\
\theta^k_j(x)&= \frac{2}{L}\sum_{\frac{(k-1)L}{2}<n<\frac{kL}{2}}(-1)^{\frac{k}{2}}\sin(\frac{2\pi n(x-j)}{L})+\frac{\sqrt{2}}{L}\sum_{n=\frac{(k-1)L}{2},\frac{kL}{2}}(-1)^{\frac{k}{2}}\sin(\frac{2\pi n(x-j)}{L}) \quad \hbox{for $k$ even}, \notag \\
% \label{equation:SOPWbasis} \\
\theta^k_j(x)&= \frac{2}{L}\sum_{\frac{(k-1)L}{2}<n<\frac{kL}{2}}(-1)^{\frac{k-1}{2}}\cos(\frac{2\pi n(x-j)}{L})+\frac{\sqrt{2}}{L}\sum_{n=\frac{(k-1)L}{2},\frac{kL}{2}}(-1)^{\frac{k-1}{2}}\cos(\frac{2\pi n(x-j)}{L}) \quad \hbox{for $k$ odd}, \notag
\end{align}}
 and in closed form:
{\footnotesize \begin{align}
\theta^1_j(x)&=\frac{1}{L}\frac{\sin(2\pi (\frac{L-1}{2})(x-j)/L)}{\sin(\pi (x-j)/L)}+\frac{\sqrt{2}}{L}\cos(\frac{2\pi(L/2)(x-j)}{L}), \notag \\
\theta^k_j(x) &=\frac{2}{L}(-1)^{\frac{k}{2}}\sin\left(\frac{2\pi (\frac{kL}{2}-\frac{L}{4})(x-j)}{L}\right)\left [ \frac{\sin(\pi(\frac{L}{2}-1)(x-j)/L)}{\sin(\pi (x-j)/L)} +\sqrt{2} \cos(\frac{\pi (x-j)}{2})\right]  \quad \hbox{for $k$ even}, \notag \\
%\label{equation:SOPWclosedForm}\\
\theta^k_j(x) &=\frac{2}{L}(-1)^{\frac{k-1}{2}}\cos\left(\frac{2\pi (\frac{kL}{2}-\frac{L}{4})(x-j)}{L}\right)\left [ \frac{\sin(\pi(\frac{L}{2}-1)(x-j)/L)}{\sin(\pi (x-j)/L)} +\sqrt{2} \cos(\frac{\pi (x-j)}{2})\right]  \quad \hbox{for $k$ odd}. \notag
\end{align}}

Equations \eqref{equation:SOPWfourierLine1}, \eqref{equation:SOPWfourierLine2} and equations \eqref{equation:fourierSOPWLine1}, \eqref{equation:fourierSOPWLine2}, \eqref{equation:fourierSOPWLine3} suggest that FFT can be used to switch between Fourier basis and SOPWs efficiently and easily. This is very important for computational purposes as it provides an efficient method to represent functions in terms of SOPWs.

Another important property of the SOPWs is that $\theta^i_j$ are the solutions to a specific variational problem. This is shown in appendix \ref{section:variation}. The final important property of $\theta^i_j$ is that for any $i$ and $j$
\[ \partial_{xx} \theta^i_j \in span \{\theta^{i}_k\}_{k=0}^{k=L-1}.\]
This can be seen by direct computation (see appendix \ref{section:firstSecondDerivative}) or using Euler-Lagrange equations for the variational problem (i.e. see appendix \ref{section:sameDepth}).

A disadvantage that SOPWs have, in comparison to Fourier basis \eqref{equation:exponentialBasis} is that SOPWs are not eigenfunctions of the derivative operator. Nevertheless, it is shown in appendix \ref{section:firstSecondDerivative} that for any $i$ and $j$,  
\[ \partial_x \theta^i_j \in span \{\theta^{i-1}_k,\theta^i_k,\theta^{i+1}_k\}_{k=0}^{L-1}.\]

%On the other hand, the great advantage that SOPWs basis have, in comparison to Fourier basis, is that shift orthogonality constraints are easy to express when we write the function in terms of SOPWs basis. Suppose $f\in L^2(\Omega)$ and 
%\begin{equation} f=\sum_{\begin{subarray}{c}i=1\\j=0,\ldots,L-1\end{subarray}}^{i=\infty} a^i_j \theta^i_j.
%\label{equation:fSOPWexpansion} \end{equation}
%Then
%\begin{equation} \int f(x)^*f(x-s) dx=\sum_{\begin{subarray}{c}i=1\\j=0,\ldots,L-1\end{subarray}}^{i=\infty} {a^i_j}^*a^i_{j+s}. \label{equation:shiftIntegral}\end{equation}
% Contrast the above expression with
%\[  \int f^*(x)f(x-s)dx=\sum_{n=-\infty}^{\infty} \hat{f}(n)^*\hat{f}(n)e^{-i2\pi s n/L}, \] 
%which arises when we expand $f$ in terms of its Fourier basis: $f(x)=\sum_{n=-\infty}^{\infty} \hat{f}(n)\phi_n(x).$ 

\section{Application to solving CPWs}\label{section:applicationCPW}
This section outlines how the projection algorithm described in section \ref{section:projectionAlgorithm} is used to compute Compressed Plain Waves (CPWs) (i.e. see \cite{Ozolins:2013CPWs}). The shift orthogonality constraints in the construction of CPWs makes their computation challenging and numerically inefficient. However, applying the projection algorithm circumvents these difficulties.

Basic compressed plane waves $\{\psi^n\}_{n=1}^\infty$ are defined by:
\begin{eqnarray}
%\left\{\begin{array}{c}
\psi^1(\x) = \displaystyle\arg\min_{\psi}\frac{1}{\mu} \int_{\Omega} |\psi(\x)| \dx + \int_{\Omega}  \psi(\x) \hat{H}_0\psi \dx  \hspace{1.5cm} \nonumber \\
           \quad \text{s.t.} \quad
            \int_{\Omega}  \psi(\x)\psi(\x - \j\w)\dx = \delta_{\j0}, \quad  \j\in\mathbb{Z}^d,  \hspace{.5cm}  \label{eqn:CPW1}
\end{eqnarray}
where $\hat{H}_0 = -\frac{1}{2}\Delta$. The higher modes can be recursively defined as:
\begin{eqnarray}
\psi^{n+1}(\x) =\displaystyle \arg\min_{\psi} \frac{1}{\mu} \int_{\Omega}  |\psi(\x)| \dx + \int_{\Omega}  \psi \hat{H}_0\psi(\x) \dx \hspace{2.2cm} \nonumber \\
       \text{s.t.}  \quad
           \left\{\begin{array}{cc}\displaystyle \int_{\Omega}  \psi(\x)\psi(\x - \j\w)\dx = \delta_{\j0}, &  \j\in\mathbb{Z}^d  \vspace{0.2cm}\\ 
           \displaystyle\int_{\Omega} \psi(\x)\psi^i(\x - \j\w)\dx = 0,  & i = 1,\cdots, n.
           \end{array}\right.  \label{eqn:CPWn}
%\end{array}\right. \\
\end{eqnarray}

To simplify our discussion, we only consider $\Omega=[0,L_1]\times [0,L_2]$, in 2D with periodic boundary conditions; the algorithms below can be straightforwardly extended to other dimensions. We use SOPWs given by equations \eqref{equation:SOPWfourierLine1} and \eqref{equation:SOPWfourierLine2} as the SOBFs used in section \ref{section:projectionAlgorithm}. In particular, we expand the given function $g$ in terms of SOPWs:
\[ g(x_1,x_2)=\sum_{i_1,i_2=1}^{N_1,N_2}\sum_{j_1=0,j_2=0}^{L_1-1,L_2-1} b_{j_1,j_2}^{i_1,i_2}\theta_{j_1}^{i_1}(x_1)\theta_{j_2}^{i_2}(x_2). \]
Operators $\SOPW$ and $\ISOPW$ are defined in the similar way to definitions \eqref{equation:definitionSOBFoperator} and \eqref{equation:definitionISOBFoperator}: 
\[ \SOPW(g)(i_1,i_2;j_1,j_2):= b_{j_1,j_2}^{i_1,i_2}, \]
and
\[ \ISOPW(b):= \sum_{i_1,i_2=1}^{N_1,N_2}\sum_{j_1=0,j_2=0}^{L_1-1,L_2-1} b_{j_1,j_2}^{i_1,i_2}\theta_{j_1}^{i_1}(x_1)\theta_{j_2}^{i_2}(x_2). \] 

Other examples of SOBFs could also be used . For this application, SOPWs were chosen mainly due to the efficiency in calculating the result of operators $\SOPW$ and $\ISOPW$ from Fourier coefficient (i.e. recall from section \ref{section:definitionSOPW} that FFT and its inverse provide an efficient procedure to switch between representation of a function in Fourier basis $\phi_n$'s and its representation in SOPWs $\theta^i_j$'s).

By introducing an auxiliary variable $u = \psi, v= \psi$, the constrained optimization problem is equivalent to the following problem:
\begin{eqnarray}
\psi^1 = \arg \min_{\psi,u} \frac{1}{\mu} \int |u(\x)| \dx  +  \int \psi \hat{H}_0\psi \dx \hspace{1.5cm} \nonumber \\
\text{s.t.} \quad  u = \psi,  v = \psi \quad \&\quad \int v(\x)v(\x - \j \w)\dx = \delta_{\j0}, \quad  j \in\mathbb{Z}^2,
\end{eqnarray}
which can be solved by an algorithm based on the Bregman iteration (i.e. see \cite{BIB:Bregman1,BIB:Bregman2, BIB:Bregman3}).

\begin{algorithm}
\label{alg:SOC_SOPW_CPW1}
\caption{Solving the first CPW using the projection algorithm}
Initialize $u^0 = v^0 = \psi^{1,0}, D^0  = B^0 = 0$.

\While{``not converged"}{
$\displaystyle \psi^{1,k}  %$ solves $(\hat{H}_0 + \lambda + r) \psi^{1,k} = \lambda( u ^{k-1} - D^{k-1} ) + r ( v^{k-1} - B^{k-1})$\;
     = \arg\min_{\psi} \int \psi \hat{H}_0 \psi \dx + \frac{\lambda}{2}\int (\psi - u^{k-1} + D^{k-1})^2\dx + \frac{r}{2}\int (\psi - v^{k-1} + B^{k-1})^2\dx$\; 
       
$\displaystyle v^{k} %= \ISOPW(Proj_{\mathcal{SSO}(N)} ( \SOPW( \psi^{1,k} + B^{k-1})))$\;
       = \arg\min_{v}\frac{r}{2}\int (\psi^{1,k} - v +B^{k-1})^2\dx$,
        \text{\rm s.t.} \, $\displaystyle \int v(\x)v(\x - \j\w)\dx = \delta_{\j0}, \quad  \j \in\mathbb{Z}^2$\;
$ \displaystyle u^k = \arg\min_{u} \frac{1}{\mu}\int  |u| \dx + \frac{\lambda}{2}\int (\psi^{1,k} - u + D^{k-1})^2 \dx $\;
          %%= \text{sgn}(\psi^{1,k}+ D^{k-1}) \max(0,|\psi^{1,k}+ D^{k-1}|-\frac{1}{\lambda\mu})$\;
            %% = \text{\rm Shrink}(\psi^{1,k}+ b^{k-1},\frac{1}{\lambda\mu})$
$D^{k} = D^{k-1} +  \psi^{1,k} - u^k$\;
$B^{k} = B^{k-1} + \psi^{1,k} -v^k$.
}
\end{algorithm}
All the above sub-optmization problems can be efficiently solved as follows:
\begin{eqnarray}
(\hat{H}_0 + \lambda + r) \psi^{1,k} = \lambda( u ^{k-1} - D^{k-1} ) + r ( v^{k-1} - B^{k-1})  \\
v^{k} = \ISOPW(Proj_{\mathcal{SSO}(N_1N_2)} ( \SOPW( \psi^{1,k} + B^{k-1}))) \\
u^k = \text{sgn}(\psi^{1,k}+ D^{k-1}) \max(0,|\psi^{1,k}+ D^{k-1}|-\frac{1}{\lambda\mu})
\end{eqnarray}

Similarly, $\psi^{n+1}$ is obtained by solving the optimization problem \eqref{eqn:CPWn} efficiently. Suppose that the first $n$ levels $\Psi^n = \{\psi^1,\cdots,\psi^n\}$ are already constructed and let ${a}^m = \SOPW(\psi^m), m = 1,\cdots,n$. In this case, the goal is to find $v^k$ satisfying
\begin{equation}
\displaystyle v^{k} = \arg\min_{v }\int (\psi^{n+1,k} - v +B^{k-1})^2 \dx,  \text{ s.t. }  
           \left\{\begin{array}{cc}\displaystyle \int  v(\x)v(\x - \j\w)\dx = \delta_{\j0}, &  \j\in\mathbb{Z}^2  \vspace{0.2cm}\\ 
           \displaystyle\int v(\x)\psi^m(\x -\j\w)\dx = 0,  & m = 1,\cdots, n.
           \end{array}\right. 
\end{equation}
Define 
\[ \mathcal{S}(\Psi^n)=span \{S(s_1,s_2){a}^m\}_{s_1=0,s_2=0,m=1}^{s_1=L_1-1,s_2=L_2-1,m=n}. \]
Using the SOPWs basis, the above problem is equivalent to solving the following problem in SOPWs frequency space: 
\begin{equation} 
Proj_{\mathcal{SSO}(N_1N_2)\cap \mathcal{S}(\Psi^n)^{\perp}} ( {b} ) := \argmin \|{b}- {v} \|_2 \quad \hbox{s.t.} \quad {v}\in \mathcal{SSO}(N_1N_2) \cap  \mathcal{S}(\Psi^n)^{\perp}.
\label{equation:SORandLinear} 
\end{equation}
%where $\SOPW(\Psi)^{\perp} =\{ v\in \mathbb{R}^{NL} ~|~ v^T S^k{a}^m=\delta_{k0}, \hbox{ for }k=0,\ldots,L-1,m = 1,\cdots,n\}$.
Vector $b$ is given and the objective is to find vector $v$ closest to $b$ that is shift orthogonal and perpendicular to ${a}^1$ to ${a}^n$.

Theorems \ref{theorem:3equivalence2D} and \ref{theorem:shiftPerpendicular2D} yield that in order to solve problem \eqref{equation:SORandLinear}, for each $j_1$ and $j_2$ one needs to find vector $z_{j_1,j_2}$ that is closest to $\B({b})(:,:;j_1,j_2)$, perpendicular to $\B({a}^m)(:,:;j_1,j_2)$ for $m=1,\ldots,n$, and lives on the unit sphere. Moreover note that, again by theorems \ref{theorem:3equivalence2D} and \ref{theorem:shiftPerpendicular2D}, $\{\B({a}^m)(:,:;j_1,j_2)\}_{m=1}^{m=n}$ form an orthonormal set of vectors for each $j_1$ and $j_2$, because elements of $\Psi^n$ are constructed such that they are shift orthogonal and orthogonal to shift span of each other. Hence, $z_{j_1,j_2}$ can be computed in two steps:
\begin{itemize}
\item $z_{j_1,j_2}=\B({b})(:,:;j_1,j_2)-\sum_{m=1}^{n}\langle \B({a}^m)(:,:;j_1,j_2), \B({b})(:,:;j_1,j_2) \rangle\B({a}^m)(:,:;j_1,j_2)$,
\item $z_{j_1,j_2}=z_{j_1,j_2}/||z_{j_1,j_2}||_2$   \qquad   (if $\|z_{j_1,j_2}\|_2\neq 0$).
\end{itemize}
In summary:

\begin{algorithm}[H]
\SetKwInOut{Input}{input}\SetKwInOut{Output}{output}
\SetKwComment{Comment}{}{}
\KwIn{${b}, {a}^1,\ldots,{a}^n$}
\KwOut{${v} = Proj_{\mathcal{SSO}(N_1N_2)\cap \mathcal{S}(\Psi^n)^{\perp}} ({b}) $}
\BlankLine
%\dontprintsemicolon
\For{$i_1=1,\ldots,N_1$ and $i_2=1,\ldots,N_2$}{
$\B(b)(i_1,i_2;~:~,~:~)=L_1L_2 \mathcal{F}_{2D}^{-1}(b(i_1,i_2;~:~,~:~))$   
}
%Set $B=W^\dagger \circulant_N(b) W_N$. \;
\For{$j_1=0,\ldots,L_1-1$ and $j_2=0,\ldots,L_2-1$}{
%Set $P_i=\sqrt{\sum_{k=1}^{N} |B_{i,(k-1)L+i}|^2}$\;
 $z_{j_1,j_2}=\B({b})(:,:;j_1,j_2)-\sum_{m=1}^{n}\langle \B({a}^m)(:,:;j_1,j_2), \B(b)(:,:;j_1,j_2)\rangle\B({a}^m)(:,:;j_1,j_2)$
\eIf {$\|z_{j_1,j_2}\|_2\neq 0$}{
	$\B(v)(:,:;j_1,j_2)=z_{j_1,j_2}/\|z_{j_1,j_2}\|_2$  
}{
	$\B(v)(:,:;j_1,j_2)$=a unit real vector orthogonal to $\{\B(a)(:,:;j_1,j_2)\}_{m=1}^{m=n}$.
}
}
\For{$i_1=1,\ldots,N_1$ and $i_2=1,\ldots,N_2$}{
$ v(i_1,i_2;~:~,~:~)=\frac{1}{L_1L_2} \mathcal{F}_{2D}(\B(v)(i_1,i_2;~:~,~:~))$  
}
%$a$ is the first row of matrix $WAW_N^\dagger$.
\caption{Projection to $\mathcal{SSO}(N_1N_2)\cap \mathcal{S}(\Psi^n)^{\perp}$}
\label{algorithm:procedureNandlinear}
\end{algorithm}

%\begin{algorithm}[H]
%\SetKwInOut{Input}{input}\SetKwInOut{Output}{output}
%\SetKwComment{Comment}{}{}
%\KwIn{${b}$}
%\KwOut{${v} = Proj_{\mathcal{SSO}(N)\cap \mathcal{S}(\Psi^n)^{\perp}} ( {b} ) $}
%\BlankLine
%%\dontprintsemicolon
%\For{$k=1,\ldots,N$}{
%$\B({b})(:,k)=F b((k-1)L+1:kL)$\;
%}
%%Set $B=W^\dagger \circulant_N(b) W_N$. \;
%\For{$i=1\ldots,L$}{
%Set $z_{i}=b_{i}-\langle b_{i},u^1_{i}\rangle u^1_{i}-\langle b_{i},u^2_{i}\rangle u^2_{i}-\cdots-\langle b_{i},u^n_{i}\rangle u^n_{i}$   \tcp*[r]{$b_i=\B({b})(i,:).$} 
%\eIf {$\|w_i\|_2\neq 0$}{
%	$\B({v})(i,:)= w_i/\|w_i\|_2$\;
%	}{
%	$\B({v})(i,:) = $ a real valued unit vector perpendicular to $\{u^1_{i}, \cdots,u^n_{i}\}$.
%}
%}
%\For{$k=1,\ldots,N$}{
%$ v((k-1)L+1:kL)=F^{-1}\B({v})(:,k)$\;
%}
%%$a$ is the first row of matrix $WAW_N^\dagger$.
%\caption{Projection to $\mathcal{SSO}(N)\cap S(\Psi^n)^{\perp}$}
%\label{algorithm:procedureNandlinear}
%\end{algorithm}

Therefore, the following algorithm is used to solve for $\psi^{n+1}$:
\begin{algorithm}
\label{alg:SOC_SOPW_CPWn}
\caption{Solving the n+1-th BCPW using the projection algorithm}
Initialize $u^0 = v^0 = \psi^{n+1,0}, D^0  = B^0 = 0$.

\While{``not converged"}{
$\displaystyle \psi^{n+1,k} $ solves $(\hat{H}_0 + \lambda + r) \psi^{n+1,k} = \lambda( u ^{k-1} - D^{k-1} ) + r ( v^{k-1} - B^{k-1})$\;

$\displaystyle v^{k} =  \ISOPW(Proj_{\mathcal{SSO}(N_1N_2)\cap \mathcal{S}(\Psi^n)^{\perp}} ( \SOPW( \psi^{n+1,k} + B^{k-1})))$\;

$ \displaystyle u^k = \text{sgn}(\psi^{1,k}+ D^{k-1}) \max(0,|\psi^{n+1,k}+ D^{k-1}|-\frac{1}{\lambda\mu})$\;
            %% = \text{\rm Shrink}(\psi^{1,k}+ b^{k-1},\frac{1}{\lambda\mu})$
$D^{k} = D^{k-1} +  \psi^{n+1,k} - u^k$\;
$B^{k} = B^{k-1} + \psi^{n+1,k} -v^k$.
}
\end{algorithm}

\begin{figure}[h]
\begin{minipage}{0.49\linewidth}
\centering
\includegraphics[width=1\linewidth]{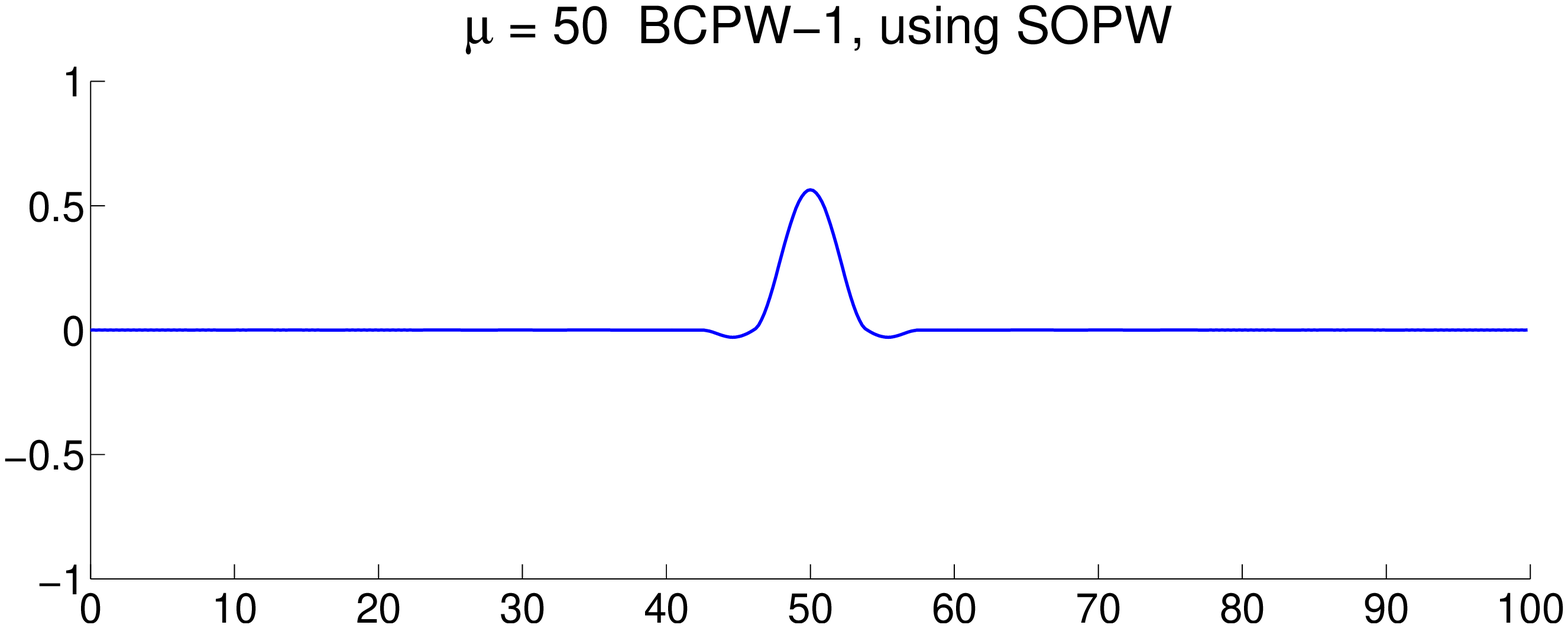}\\
\end{minipage}\hfill
\begin{minipage}{0.49\linewidth}
\centering
\includegraphics[width=1\linewidth]{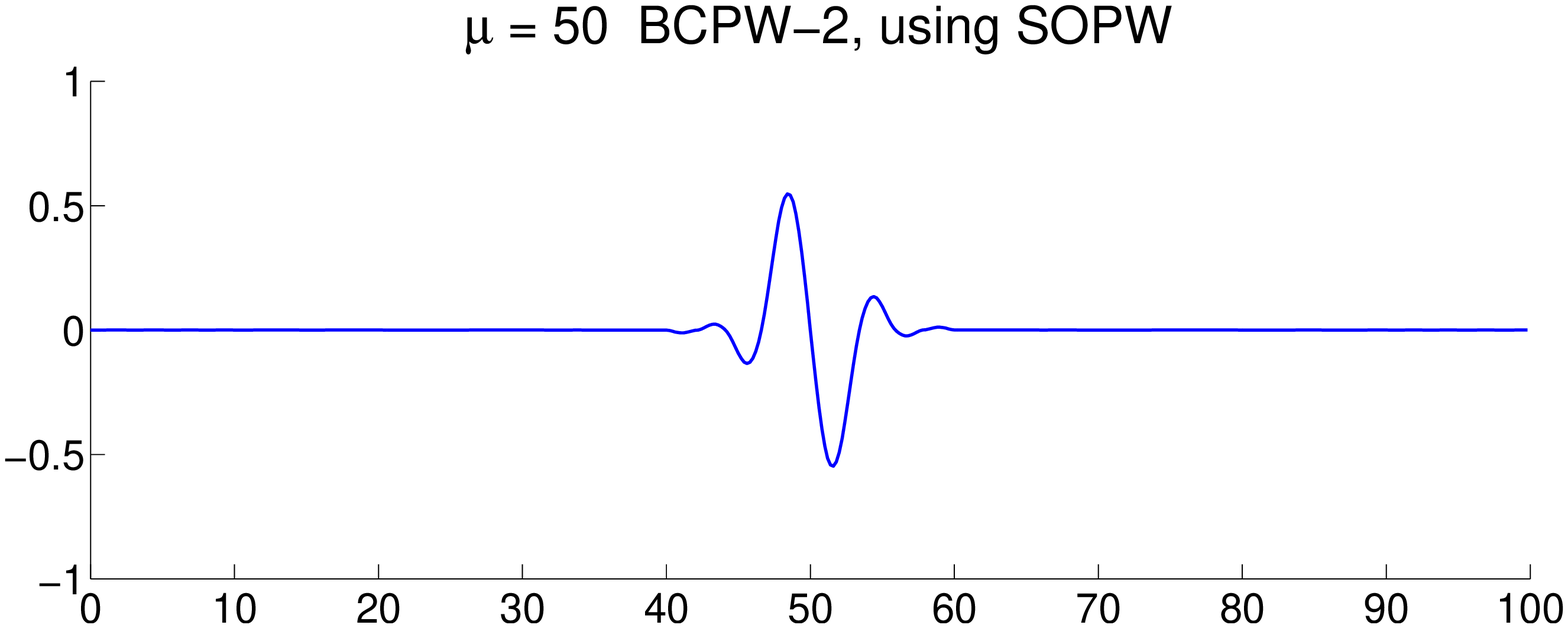}\\
\end{minipage}\hfill\\
\begin{minipage}{0.49\linewidth}
\centering
\includegraphics[width=1\linewidth]{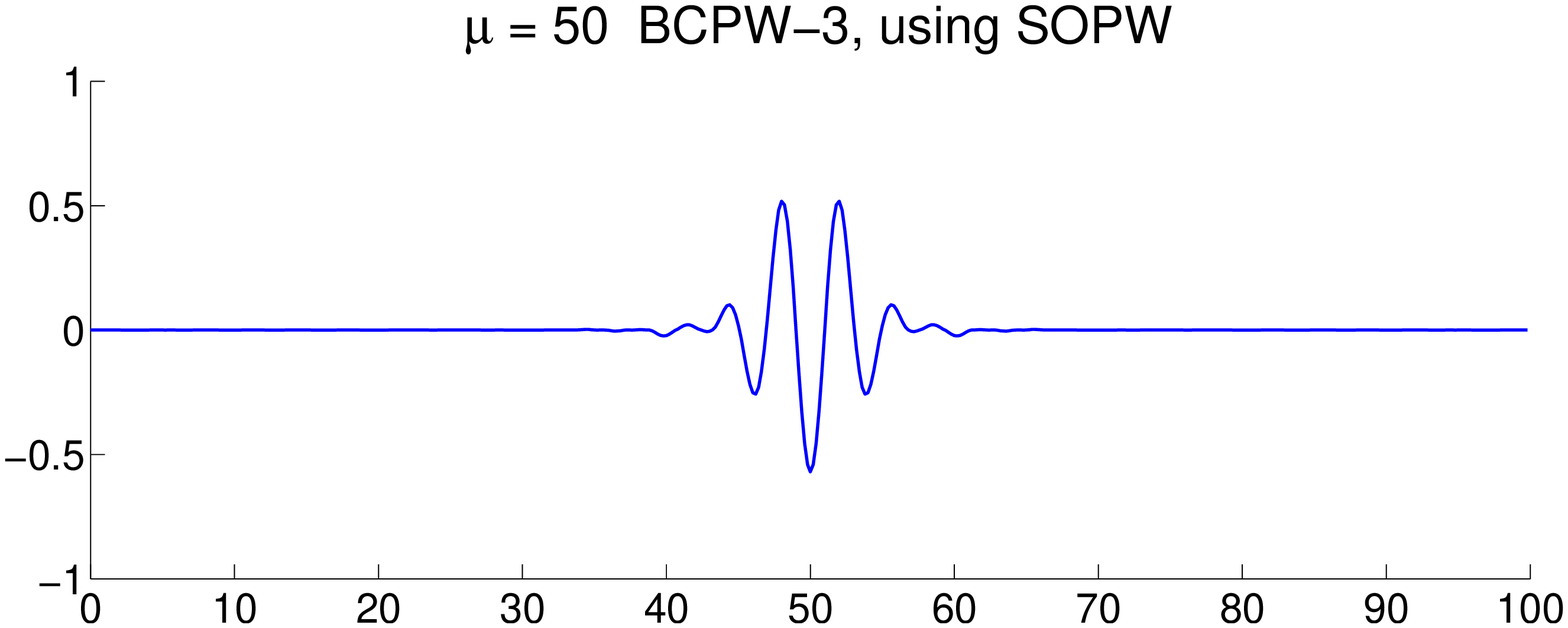}\\
\end{minipage}\hfill
\begin{minipage}{0.49\linewidth}
\centering
\includegraphics[width=1\linewidth]{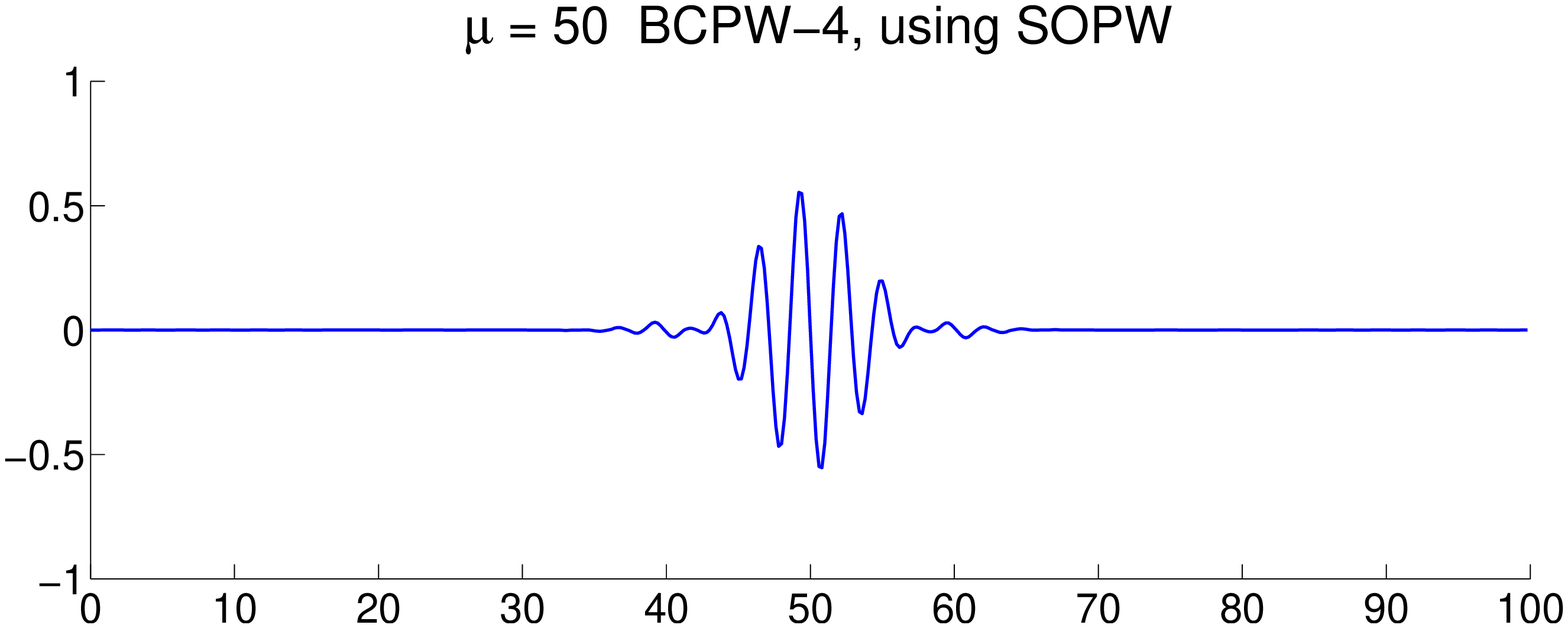}\\
\end{minipage}\hfill
\caption{The first 4 one-dimensional BCPWs obtained by using the new procedure.}
\label{fig:BCPWs}
\end{figure}

Figure \ref{fig:BCPWs} plots the first four BCPWs in 1D using the proposed algorithms. These results are very consistent with the results in \cite{Ozolins:2013CPWs}. Table \ref{tab:EfficiencyComparisons}, highlights the computational speed gained by using the new procedure outlined in this section.

\begin{table}[h]
\begin{center}
\begin{tabular}{|c||c|c|c|c|c||c|c|c|c|c|}
\hline
\multirow{2}{1.1cm}{$\#$ of points}&   \multicolumn{5}{|c||}{Algorithm proposed in \cite{Ozolins:2013CPWs}} &  \multicolumn{5}{|c|}{The new procedure}  \\
 \cline{2-11}
  & $\psi^1$ & $\psi^2$ & $\psi^3$ & $\psi^4$ & total &  $\psi^1$ & $\psi^2$ & $\psi^3$ & $\psi^4$ & total \\
  \hline
500 &  23.05 & 29.39 & 21.40 & 9.42 &  83.26 & 0.30 & 1.47 & 0.44 & 0.28 &  2.49    \\
 \hline
1000   & 53.11 & 116.30 & 81.28 & 29.99 & 280.68   & 0.71 & 3.48 & 1.23 & 0.66 &  6.08    \\
 \hline
\end{tabular} \vspace{-0.5cm}
\end{center}
\caption{CPU time consumption (seconds) for computing the first 4 one-dimensional BCPWs using SOPWs and FFT with the same accuracy.}
\label{tab:EfficiencyComparisons}
\end{table}

\section{Conclusion}\label{section:conclusion}
This paper presents a fast algorithm for finding a closest shift orthogonal function to a given function. The algorithm can be easily implemented using FFT and has computational complexity bounded by $M\log(M)$, where $M$ is the number of coefficients used to store the input function. 

The algorithm described here is very useful for problems with shift orthogonality constraints. As an example, the application of the algorithm in computing Compressed Plain Waves (CPWs) is shown.  

\pagebreak

\appendix

\section{Variational Origin of SOPWs}\label{section:variation}

Here, it is shown that  in $1D$, if there is no $L^1$ term in the definition of CPWs (i.e. $\mu=\infty$) then they are essentially the same as SOPWs given by \eqref{equation:SOPWfourierLine1} and \eqref{equation:SOPWfourierLine2}. Let $\Omega=[0,L]$ (where $L$ is even) and by scaling assume that $w=1$. For any function $\psi$ define
\begin{equation}
\mathcal{J}_\infty(\psi):=\int_\Omega \psi\hat{H}_0\psi \dx .
\end{equation}
As before, $\hat{H}_0=-\frac{1}{2}\partial_{xx}$. It is clear that $\hat{H}_0$ has eigenfunctions $\phi_n(x)=\frac{1}{\sqrt{L}}e^{i2\pi n x/L}$ with corresponding eigenvalue $\lambda_n=2(\pi n/L)^2$, $n=0,\pm1,\pm 2,\ldots$.
Note that BCPWs $\{ \theta^i\}_{i=1}^{i=\infty}$ when there is no $L^1$ term (i.e. $\mu=\infty$) are defined in the following way:
%\begin{equation}
%\begin{cases}
%\theta^1=\argmin_{\theta} \mathcal{J}_{\infty}(\theta) \quad \hbox{ s.t. } \quad \int \theta(x)\theta(x-j)dx=\delta_{j0} \\
%\theta^k=\argmin_{\theta} \mathcal{J}_{\infty}(\theta)  \quad \hbox{ s.t. } \quad \int \theta(x)\theta(x-j)dx=\delta_{j0}  \quad\&\quad \int \theta(x)\theta^i(x-j)dx=0
%\end{cases}
%\label{equation:xiProblem}\end{equation} 
\begin{align}
\psi^1&=\argmin_{\psi} \mathcal{J}_{\infty}(\psi) \quad \hbox{ s.t. } \quad \int \psi(x)\psi(x-j)dx=\delta_{j0}, \notag \\
\psi^k&=\argmin_{\psi} \mathcal{J}_{\infty}(\psi)  \quad \hbox{ s.t. } \quad\begin{cases} \int \psi(x)\psi(x-j)dx=\delta_{j0}, \\
\int \psi(x)\psi^i(x-j)dx=0 \quad \hbox{for }i=1,\ldots, k-1.
\end{cases}
\label{equation:xiProblem}\end{align}

The main result of this section is the following theorem, which implies $\{ \psi^i\}_{i=1}^{i=\infty}$ are indeed SOPWs  $\{ \theta^i\}_{i=1}^{i=\infty}$:
\begin{theorem}
One set of solutions to problem \eqref{equation:xiProblem} are
\[  \theta^1(x)=\frac{1}{\sqrt{L}}\sum_{|n|<\frac{L}{2}}\phi_n(x)+\frac{1}{\sqrt{2L}}\sum_{|n|=\frac{L}{2}}\phi_{n}(x) ,  \]
and, for $k>1$ 
\[ \theta^k(x)=\frac{1}{\sqrt{L}}\sum_{\frac{(k-1)L}{2}<|n|<\frac{kL}{2}}(\sgn(n)i)^{k-1}\phi_n(x)+\frac{1}{\sqrt{2L}}\sum_{|n|=\frac{(k-1)L}{2},\frac{kL}{2}}(\sgn(n)i)^{k-1}\phi_n(x)  \]
\label{theorem:xiFormula}\end{theorem}
\proof Since $\phi_n$'s form complete set, for any function $\theta(x)\in L^2$ we can write
\begin{equation} \theta(x)=\sum_{n=-\infty}^{\infty} a(n)\phi_n(x).\label{equation:eigenvecExpansion}\end{equation}
%Since we are interested only in real functions, we must have 
%\begin{equation}  a(-n)=a^*(n), \qquad \hbox{for } n=1,2,\ldots \label{equation:realness}\end{equation}
%In view of \eqref{equation:eigenvecExpansion} and \eqref{equation:realness}, 
Consequently, 
\begin{equation} \mathcal{J}_{\infty}(\theta)=\sum_{n=-\infty}^{\infty}|a(n)|^2 \lambda_n=\sum_{n=1}^{\infty}(|a(n)|^2+|a(-n)|^2)\lambda_n, \label{equation:objectivePre}\end{equation}
and shift orthogonality constraints yield that
\begin{align} 
&\delta_{j0}=\int \theta^*(x)\theta(x-j)dx=\sum_{n=-\infty}^{\infty} |a(n)|^2e^{-i2\pi jn/L} \notag \\
&=|a(0)|^2+\sum_{n=1}^{\infty}(|a(n)|^2+|a(-n)|^2)\cos(2\pi jn/L)-i \sum_{n=1}^{\infty}(|a(n)|^2-|a(-n)|^2)\sin(2\pi jn/L). 
\label{equation:constraintsPre} \end{align}
Therefore, for $\theta(x)$ to be feasible (i.e satisfy shift orthogonality constraints), it must be the case that $|a(n)|=|a(-n)|$ for all $n\geq 1$. Note that changing the phase value of $a(n)$ and $a(-n)$ does not change the value of the objective function \eqref{equation:objectivePre}. Thus, the phase factor of $a(n)$ and $a(-n)$ can be chosen in such a way that 
\begin{equation}  a(-n)=a^*(n), \qquad \hbox{for } n=1,2,\ldots. \label{equation:realness}\end{equation}
The above conditions guarantee that $\theta(x)$ is real valued. Hence, there always exist a real valued minimizers for variational problem \eqref{equation:xiProblem}. In view of \eqref{equation:realness}, equations \eqref{equation:objectivePre} and \eqref{equation:constraintsPre} can be re-written: objective function becomes
\begin{equation} \mathcal{J}_{\infty}(\theta)=2\sum_{n=1}^{\infty}|a(n)|^2 \lambda_n, \label{equation:objective}\end{equation}
and shift orthogonality constraints yield that for $j=0,1,\ldots,L-1,$
\begin{equation} \delta_{j0}=|a(0)|^2+2\sum_{n=1}^{\infty}|a(n)|^2\cos(2\pi jn/L). \label{equation:constraints} \end{equation}
Lets first find $\theta^1$. To that end, the goal is to find $\{a(n)\}_{n=0}^\infty$ that minimizes \eqref{equation:objective} and satisfies \eqref{equation:constraints} for $j=0,\ldots,L-1$. Set 
\begin{equation}
 \begin{cases} c(0)=|a(0)|^2 \\
c(n)=2|a(n)|^2 \qquad \hbox{for } n=1,2,\dots.
\end{cases}
\label{equation:ca}\end{equation}

Let $M$ be the $L\times L$ matrix whose $(j,n)$-th entry is $\cos(2\pi jn/L)$ for $j,n=0,\ldots,L-1$. Form infinite dimensional matrix $A$ by concatenating infinitely many copies of $M$ side by side, that is 
\[ A=[M|M|\cdots]. \]
Therefore, to find $\theta^1$, one needs to solve the following optimization problem:\
\begin{equation} \argmin_c  \lambda^Tc \qquad \hbox{s.t.} \qquad Ac=b, c\geq 0, \label{equation:minimizationProblem}\end{equation}
where 
\[  \lambda^T=[\lambda_0,\lambda_1,\cdots], \quad c^T=[c(0),c(1),\cdots], \quad \hbox{and} \quad b^T=[1,\underbrace{0,\cdots,0}_{L-1}]. \] 

Observe that matrix $M$ is not invertible. However, it can be partitioned in the following way:
Let $\vec{1}=M(:,1)$, $M_L=M(:,2:L/2-1)$, $\vec{e}=M(:,L/2)$ and $M_R=M(:,L/2+1:L-1)$. Then 
\[ A=[\vec{1}|M_L|\vec{e}|M_R|\vec{1}|M_L|\vec{e}|M_R|\cdots]. \]

The strategy is to guess the solution to problem \eqref{equation:minimizationProblem} and then verify (i.e. using the dual formulation of \eqref{equation:minimizationProblem}) that it is indeed the optimal solution. For this purpose we first prove the following three lemmas:
\begin{lemma} The $L\times (L/2-1)$ matrices $[\vec{1}|M_L|\vec{e}]$ and $[\vec{e}|M_R|\vec{1}]$ have full rank.
\label{lemma:invertible}
\end{lemma}
\proof First observe that matrix $[\vec{e}|M_R|\vec{1}]$ is formed from matrix $[\vec{1}|M_L|\vec{e}]$ if the columns are arranged in the opposite order. Thus it suffices to only show $[\vec{1}|M_L|\vec{e}]$ has full rank. For contrary assume the opposite that $[\vec{1}|M_L|\vec{e}]$ is not full rank, then there exist nontrivial set of constants $\{k_n\}_{n=0}^{n=L/2}$ such that 
\[ \sum_{n=0}^{L/2}k_n\cos(2\pi jn/L)=0 \qquad \hbox{for } j=0,\ldots,L-1. \]
The above system of equations implies that 
\begin{equation}  \sum_{n=0}^{L-1}k'_n e^{-i2\pi jn/L}=0 \qquad \hbox{for } j=0,\ldots,L-1, \label{equation:colsDependance} \end{equation}
where 
\[ k'_n=\begin{cases}
k_n \qquad &\hbox{for } n=0 \hbox{ and } L/2,\\
k_n/2  \qquad &\hbox{for } 0<n<L/2, \\
k_{L-n}/2  \qquad &\hbox{for } L/2<n\leq L-1. 
\end{cases} \] 
However, system of equations \eqref{equation:colsDependance} implies that the columns of the $L\times L$ Discrete Fourier Transform matrix are linearly dependent; which contradicts invertibility of the DFT matrix.\qed

\begin{lemma} For $k\geq 1$ and $j=0,\ldots,L-1$:
\[ \frac{1}{L}\cos(\pi j(k-1))+\frac{2}{L}\sum_{\frac{(k-1)L}{2}<n<\frac{kL}{2}}\cos(2\pi jn/L)+\frac{1}{L}\cos(\pi jk)=\delta_{j0}. \]
\label{lemma:sumCos}\end{lemma}
\proof If $j=0$ the result is clear. For $j=1,\ldots,L-1$:
\begin{align}
~& \frac{1}{L}\cos(\pi j(k-1))+\frac{2}{L}\sum_{\frac{(k-1)L}{2}<n<\frac{kL}{2}}\cos(2\pi jn/L)+\frac{1}{L}\cos(\pi jk)= \notag\\
=&  \frac{1}{L}\cos(\pi j(k-1))+\frac{2}{L}\sum_{\frac{(k-1)L}{2}<n<\frac{kL}{2}}\frac{1}{2}\left[\cos(2\pi jn/L)+\cos(2\pi j(kL-n)/L)\right]+\frac{1}{L}\cos(\pi jk)= \notag \\
=& \frac{1}{L}\sum_{\frac{(k-1)L}{2}\leq n<\frac{(k+1)L}{2}}\cos(2\pi jn/L) = \frac{1}{L} Re\left\{\sum_{n=(k-1)L/2}^{(k+1)L/2-1}e^{i2\pi jn/L}\right \} = \notag \\
=& \frac{1}{L} Re\left\{ e^{i\pi j(k-1)}(\frac{1-e^{i2\pi jL/L}}{1-e^{i2\pi j/L}})\right \}=0. \notag
\end{align}
The result follows.\qed
\begin{lemma}
Assume that there exist vector $y\in \mathbb{R}^L$ and infinite dimensional vector $s$ such that 
\[ A^Ty+s=\lambda, \qquad s\geq 0. \]
If $c$ is feasible for problem \eqref{equation:minimizationProblem}, then 
%\[ \lambda^Tc-b^Ty=s^Tc\geq 0. \]
\[  \lambda^Tc \geq y^Tb.\]
Furtheremore, if additionally $s$ and $c$ satisfy complementary slackness property $s^Tc=0$, then $c$ is the solution to problem \eqref{equation:minimizationProblem}.
\label{lemma:optimality}\end{lemma}
\proof Observe that for any feasible $c$ in problem \eqref{equation:minimizationProblem},
\[ \lambda^Tc=(A^Ty+s)^Tc=y^TAc+s^Tc=y^Tb+s^Tc\geq y^Tb,\]
where the last inequality is from nonnegativity of $s$ and $c$. Hence, the minimum value of the objective function in problem \eqref{equation:minimizationProblem} is $y^Tb$. Moreover, the minimum is achieved if $c$ is feasible and satisfies $s^Tc=0$. \qed

Now returning to optimization problem \eqref{equation:minimizationProblem}, set
\[ c^T=[\frac{1}{L},\underbrace{\frac{2}{L},\ldots,\frac{2}{L}}_{L/2-1},\frac{1}{L},0,0,\ldots],  \]
and find $y$ that satisfies 
\[ [\vec{1}|M_L|\vec{e}]^T y=[\lambda_0,\ldots,\lambda_{L/2}]^T. \]
Such $y$ exist because by lemma \ref{lemma:invertible} matrix $[\vec{1}|M_L|\vec{e}]$ is full rank. Finally, set
\[ s=\lambda-A^Ty. \]
Lemma \ref{lemma:sumCos} implies that $c$ is feasible for problem \eqref{equation:minimizationProblem}. Moreover, it is straightforward to verify that
\begin{align*} s^T=[\underbrace{0,\ldots,0}_{L/2+1},&\lambda_{\frac{L}{2}+1}-\lambda_{
\frac{L}{2}-1},\lambda_{\frac{L}{2}+2}-\lambda_{
\frac{L}{2}-2} ,\ldots, \lambda_{L}-\lambda_{0},\lambda_{L+1}-\lambda_1,\ldots \\
&\ldots ,\lambda_{\frac{3L}{2}+1}-\lambda_{
\frac{L}{2}-1},\lambda_{\frac{3L}{2}+2}-\lambda_{
\frac{L}{2}-2},\ldots, \lambda_{2L}-\lambda_{0},\lambda_{2L+1}-\lambda_1,\ldots] \geq 0. \end{align*}
%\[ s^T=[\underbrace{0,\ldots,0}_{L/2+1},\lambda_{\frac{L}{2}+1}-\lambda_{
%\frac{L}{2}-1},\ldots, \lambda_{L}-\lambda_{0},\lambda_{L+1}-\lambda_1,\ldots,\lambda_{\frac{3L}{2}+1}-\lambda_{
%\frac{L}{2}-1},\ldots, \lambda_{2L}-\lambda_{0},\lambda_{2L+1}-\lambda_1,\ldots] \geq 0. \]
Thus, by lemma \ref{lemma:optimality}, $c$ is the solution of problem \eqref{equation:minimizationProblem}. Hence, in view of \eqref{equation:ca},
\[ \begin{cases}
|a(n)|=1/\sqrt{L} \qquad &\hbox{for } n=\pm 1,\ldots,\pm (\frac{L}{2}-1), \\
|a(\pm \frac{L}{2})|=1/\sqrt{2L} \qquad &~\\
|a(n)|=0 \qquad &\hbox{otherwise}.
\end{cases}\]
Note that any phase values for $a(n)$ as long as \eqref{equation:realness} holds is acceptable. If all the phase factors are set to equal to 1, then
\[  \theta^1(x)=\frac{1}{\sqrt{L}}\sum_{|n|<\frac{L}{2}}\phi_n(x)+\frac{1}{\sqrt{2L}}\sum_{|n|=\frac{L}{2}}\phi_{n}(x).  \]
Next, lets find $\theta^2$. For $j=0,\ldots,L-1$, define 
\[ \theta^1_j(x):=\theta^1(x-j). \]
Observe that from relationships \eqref{equation:fourierSOPWLine1} and \eqref{equation:fourierSOPWLine2} (i.e. with $k=1$),
%\begin{equation}  \{\phi_{-\frac{L}{2}+1},\ldots,\phi_0,\ldots,\phi_{\frac{L}{2}-1}\} \subset span\{\theta^1_j\}_{j=0}^{L-1}, \label{equation:spanSpace}\end{equation}
\[  \{\phi_n\}_{n=-L/2+1}^{n=L/2-1} \subset span\{\theta^1_j\}_{j=0}^{j=L-1}.\]
In particular, orthogonality to previous CPWs constraints in problem \eqref{equation:xiProblem}, imply that if $\theta^2$ is expanded in the form \eqref{equation:eigenvecExpansion}, it is necessary (but not sufficient) that $a(0)=\cdots=a(\pm(\frac{L}{2}-1))=0$. 
Let $\theta^*$ be the solution of problem 
\begin{equation} \min_{\theta} \mathcal{J}_{\infty}(\theta) \quad \hbox{ s.t. } \quad \int \theta(x)\theta(x-j)dx=\delta_{j0}, \label{equation:auxillaryXi}\end{equation}
with an additional constraint that if $\theta^*$ is expanded in the form \eqref{equation:eigenvecExpansion}, then $a(n)=0$ for $|n|<L/2$. 
Using the same arguments as before, one concludes that to find  a candidate for $\theta^*$ it is required to solve an optimization problem similar to \eqref{equation:minimizationProblem}; however, this time 
\[  \lambda^T=[\lambda_{\frac{L}{2}},\lambda_{\frac{L}{2}+1},\cdots], \quad c^T=[c({\frac{L}{2}}),c({\frac{L}{2}}+1),\cdots], \quad \hbox{and} \quad M=[\vec{e}|M_R|\vec{1}|M_L].  \] 
Repeating the same line of logic as before, the optimal solution is still 
\[ c^T=[\frac{1}{L},\underbrace{\frac{2}{L},\ldots,\frac{2}{L}}_{L/2-1},\frac{1}{L},0,0,\ldots].  \]
Hence, for $\theta^*$,

\[ \begin{cases}
|a(\pm \frac{L}{2})|=|a(\pm L)|=1/\sqrt{2L} \qquad &~ \\
|a(n)|=1/\sqrt{L} \qquad &\hbox{for } n=\pm  (\frac{L}{2}+1),\ldots,\pm (L-1), \\
|a(n)|=0 \qquad &\hbox{otherwise}.
\end{cases}\]

Now observe that $\theta^*$ is not necessarily the same as $\theta^2$, as $\theta^2$ satisfies stricter constraints (i.e. orthogonality to $\{\theta^1_j\}_{j=0}^{j=L-1}$) than $\theta^*$. Nevertheless, there is a particular choice of phases for $a(n)$'s for which 
\[ \theta^{*}(x) =\frac{1}{\sqrt{L}}\sum_{\frac{L}{2}<|n|<L}(\sgn(n)i)\phi_n(x)+\frac{1}{\sqrt{2L}}\sum_{|n|=\frac{L}{2},L}(\sgn(n)i)\phi_n(x).  \]
It is easy to verify that the above function is indeed orthogonal to set $\{\theta^1_j\}_{j=0}^{j=L-1}$. Therefore,
\[ \theta^{2}(x) =\frac{1}{\sqrt{L}}\sum_{\frac{L}{2}<|n|<L}(\sgn(n)i)\phi_n(x)+\frac{1}{\sqrt{2L}}\sum_{|n|=\frac{L}{2},L}(\sgn(n)i)\phi_n(x).  \]

Continue the above procedure to find the subsequent BCPWs: for example if $\theta^3$ is expanded in the form \eqref{equation:eigenvecExpansion}, it is necessary (but not sufficient) that $a(n)=0$, for $|n|<L$. For from relationships \eqref{equation:fourierSOPWLine1}, \eqref{equation:fourierSOPWLine2} (i.e. with $k=1,2$) and \eqref{equation:fourierSOPWLine3} (i.e. with $k=1$):
\[ \{\phi_n\}_{n=-(L-1)}^{n=L-1} \subset span\{\theta^1_j,\theta^2_j\}_{j=0}^{j=L-1}, \]
and from the orthogonality to previous CPWs constraints in problem \eqref{equation:xiProblem}.

Let $\theta^*$ be the solution of problem \eqref{equation:auxillaryXi} with an additional constraint that if it is expanded in the form \eqref{equation:eigenvecExpansion}, then $a(n)=0$ for $|n|<L$. We conclude that for $\theta^*$,
\[ \begin{cases}
|a(\pm L)|=|a(\pm \frac{3L}{2})|=1/\sqrt{2L} \qquad &~ \\
|a(n)|=1/\sqrt{L} \qquad &\hbox{for } n=\pm  (L+1),\ldots,\pm (\frac{3L}{2}-1), \\
|a(n)|=0 \qquad &\hbox{otherwise}.
\end{cases}\]
Note that $\theta^*$ is already orthogonal to the space spanned by all shifts of function $\theta^1$. To make $\theta^*$ orthogonal to the space spanned by all shifts of function $\theta^2$ (and therefore, derive a formula for $\theta^3$), a particular choice of phase factors for $a(n)$'s are chosen. Consequently, 
\[ \theta^{3}(x) = \frac{1}{\sqrt{L}}\sum_{L<|n|<\frac{3L}{2}}(-1)\phi_n(x)+\frac{1}{\sqrt{2L}}\sum_{|n|=L,\frac{3L}{2}}(-1)\phi_n(x). \]     

\section{Laplacian of the SOPWs}\label{section:sameDepth}
Here, it is shown that for any set of solutions to the variational problem \eqref{equation:xiProblem}:
\[ \partial_{xx} \theta^k_j \in span \{\theta^{k}_\ell\}_{\ell=0}^{\ell=L-1},\]
where $\theta^k_j(x):=\theta^k(x-j)$. Observe that it suffices to show that 
\[ \partial_{xx} \theta^k \in span \{\theta^{k}_\ell\}_{\ell=0}^{\ell=L-1}.\]
From the theory of variational calculus with constraints (i.e. see for example \cite[Chapter 8]{BIB:evans}) at the $k$-th step (i.e. when $\theta^i_j$ for $i=1,\ldots,k-1$ are already determined), if $\theta^k$ is the solution to the variational problem \eqref{equation:xiProblem}, then it is the weak solution of the Euler-Lagrange equation
\begin{equation}
\Delta \theta^k=\sum_{i=1}^{k} \sum_{j=0}^{L-1}\lambda_j^i \theta_j^i,
\label{equation:EL}\end{equation}
where constants $\lambda^i_j$ are the Lagrange multipliers corresponding to the orthonormality constraints: 
\[ \int \theta^k(x) \theta^k_j(x) \d x=\delta_{j0} \quad \hbox{ and } \int \theta^k(x) \theta^i_j(x) \d x=0 \quad \hbox{ for } i=1,\ldots,k-1. \]
It remains to show that $\lambda^i_j=0$ for all $i<k$. Fix $n<k$ and $\ell \in \{0,\ldots,L-1\}$. Multiply both sides of \eqref{equation:EL} by $\theta^n_\ell$, integrate over the domain $[0,L]$, and use orthonormalities of $\{\theta^i_j\}_{i=1,j=0}^{i=k,j=L-1}$  and integration by parts to conclude that 
\begin{equation} \lambda^n_\ell=\int \Delta \theta^k(x) \theta^n_\ell(x) \d x=\int \theta^k(x)\Delta\theta^n_\ell(x) \d x=\int \theta^k_{L-\ell}(x)\Delta\theta^n(x) \d x. \label{equation:integrationByParts}\end{equation}

Next, observe that $\theta^n$ must satisfy a similar equation to \eqref{equation:EL}; that is,  
\[ \Delta \theta^n=\sum_{i=1}^{n} \sum_{j=0}^{L-1}\gamma_j^i \theta_j^i. \]
From definition of $\theta^k$ in \eqref{equation:xiProblem} and because $n<k$, one concludes that $\theta^k_{L-\ell}$ is orthogonal to $\{\theta^i_j\}_{i=1,j=0}^{i=n,j=L-1}$. Therefore, multiplying the above equation by $\theta^k_{L-\ell}$ and integrating over the domain, yields that
\[ \int \theta^k_{L-\ell}\Delta \theta^n(x) \d x=0. \]
The above equation and equation \eqref{equation:integrationByParts} imply that $\lambda^n_\ell=0$ as was to be shown.

\section{First and Second Derivative of SOPWs}\label{section:firstSecondDerivative}
Here, formulas for the first and second derivatives of SOPWs, defined by \eqref{equation:SOPWfourierLine1} and \eqref{equation:SOPWfourierLine2}, are presented. These results are important in determining the matrix elements of the derivative and Laplacian operator when SOPWs basis are used.

The first derivative of of SOPWs are given by the following theorem:
\begin{theorem}[First Derivatives]\label{theorem:firstDerivative}
For $k=1,2,\ldots$ and $\ell=0,\ldots,L-1$,
\[ \partial_x \theta^{k}_\ell= \frac{\pi}{L}\left[ -(k-1)\sum_{j} (-1)^{(k-1)(j-\ell)}\theta^{k-1}_j+\sum_{j} a(j-\ell) \theta^{k}_j+k\sum_{j} (-1)^{k(j-\ell)}\theta^{k+1}_j \right], \]
where
\[ a(j-\ell)=\begin{cases}
0 \qquad &\hbox{if } j-\ell=0, \\
(-1)^k(2k-1)\cot(\pi(j-\ell)/L)   &\hbox{if $j-\ell$ is odd}, \\ 
\cot(\pi (j-\ell)/L) &\hbox{otherwise}, 
\end{cases} \]
and dummy variable $j$ takes its values values from $\{0,1,\ldots,L-1\}$.
\end{theorem}

The second derivative of SOPWs are given by the following theorem:
\begin{theorem}[Second Derivatives]\label{theorem:secondDerivative}
For $k=1,2,\ldots$ and $\ell=0,\ldots,L-1$,
\[ \partial_{xx}\theta^{k}_\ell =\frac{-\pi^2}{L^2}\sum_j b(j-\ell) \theta^{k}_j\]
where
\[ b(j-\ell)=\begin{cases}
(k^2-k+1/3)L^2+2/3 \qquad &\hbox{if } j-\ell=0, \\
(-1)^{k}(4k-2)\csc^2(\pi(j-\ell)/L)   &\hbox{if $j-\ell$ is odd}, \\ 
2 \csc^2(\pi(j-\ell)/L)  &\hbox{otherwise}, 
\end{cases} \]
and dummy variable $j$ takes its values from $\{0,1,\ldots,L-1\}$.
\end{theorem}

Recall that  $\omega_j=e^{i2\pi j/L}$. The following lemma is essential in the proof of the above theorems:
\begin{lemma}\label{lemma:sumOmega}
For positive integer $k$, even $L$ and $j=0,\ldots,L-1$:
\[  \sum_{\frac{(k-1)L}{2}<|n|<\frac{kL}{2}} n\omega_j^n=
\begin{cases}
0 \qquad &\hbox{if } j=0, \\
\frac{-i}{2}L (-1)^{k}(2k-1) \cot(\pi j/L)   &\hbox{if $j$ is odd}, \\ 
\frac{-i}{2} L \cot(\pi j/L)  &\hbox{otherwise}, 
\end{cases} \]
and
\[  \sum_{\frac{(k-1)L}{2}<|n|<\frac{kL}{2}} n^2\omega_j^n=
\begin{cases}
\frac{1}{12}(L-2)L[(3k^2-3k+1)L-1] \qquad &\hbox{if } j=0, \\
\frac{L}{2}(-1)^k(2k-1)\csc^2(\pi j/L)-(-1)^k(2k-1)\frac{L^2}{4}   &\hbox{if $j$ is odd}, \\ 
\frac{L}{2}\csc^2(\pi j/L)-(k^2+(k-1)^2)\frac{L^2}{4}  &\hbox{otherwise}, 
\end{cases} \]
\end{lemma}   
\proof  It is easy to verify the case $j=0$, so we assume $j\neq 0$. Let 
\[ c(x)=\sum_{\frac{(k-1)L}{2}<n<\frac{kL}{2}} \cos(nx). \]
Then 
\[ \sum_{\frac{(k-1)L}{2}<|n|<\frac{kL}{2}} n\omega_j^n=\sum_{\frac{(k-1)L}{2}<n<\frac{kL}{2}} n(\omega_j^n-\bar{\omega}_j^n)=2i \sum_{\frac{(k-1)L}{2}<n<\frac{kL}{2}} n \sin(\frac{2\pi nj}{L})=(-2i)c'(\frac{2\pi j}{L}), \]
and 
\[ \sum_{\frac{(k-1)L}{2}<|n|<\frac{kL}{2}} n^2\omega_j^n=\sum_{\frac{(k-1)L}{2}<n<\frac{kL}{2}} n^2(\omega_j^n+\bar{\omega}_j^n)=2 \sum_{\frac{(k-1)L}{2}<n<\frac{kL}{2}} n^2\cos(\frac{2\pi nj}{L})=(-2)c''(\frac{2\pi j}{L}). \]
It is well known (i.e. for example see \cite[page 290]{BIB:Hull}) that
\[ c(x)=\frac{\sin((\frac{L}{2}-1)x/2)}{\sin(x/2)}\cos\left(\frac{(2k-1)}{4}Lx\right). \]
The rest of the proof follows from straightforward but tedious calculations: one finds close formulas for $c'(x)$ and $c''(x)$, substitutes $x=2\pi j/L$ and simplifies. In particular, table \ref{table:simplification} is helpful in simplifying. \qed

\begin{table}[h]
\begin{center}
\begin{tabular}{|c||c|c|c|c|}
\hline
$j$ \hbox{ mod 4} &  $\cos\left(\frac{\pi j (2k-1)}{2}\right)$ & $\sin\left(\frac{\pi j (2k-1)}{2}\right)$ & $\cos\left(\frac{\pi j (L/2-1)}{L}\right)$ &   $\sin\left(\frac{\pi j (L/2-1)}{L}\right)$   \\
\hline
\hline
0  &   1 & 0  & $\cos(\pi j/L)$  & $-\sin(\pi j/L)$  \\
\hline
1  &  0  & $(-1)^{k+1}$ & $\sin(\pi j/L)$  & $\cos(\pi j/L)$ \\
\hline
2  &  $-1$  & 0  & $-\cos(\pi j/L)$  & $\sin(\pi j/L)$ \\
\hline
3  & 0    & $(-1)^k$  & $-\sin(\pi j/L)$  & $-\cos(\pi j/L)$ \\
\hline
\end{tabular} \vspace{-0.5cm}
\end{center}
\caption{Trigonometry identities for integers $k$, $j$ and even positive number $L$.}
\label{table:simplification}
\end{table}
The proof of theorems \ref{theorem:firstDerivative} and \ref{theorem:secondDerivative} are very similar. The idea of the proof is simple: write SOPWs basis in terms of Fourier basis using formulas \eqref{equation:SOPWfourierLine1} and \eqref{equation:SOPWfourierLine2}, take appropriate number of derivatives, and then use formulas \eqref{equation:fourierSOPWLine2} and \eqref{equation:fourierSOPWLine3} to write back the result in terms of the SOPWs basis.

\textbf{Proof of theorem \ref{theorem:firstDerivative}:} First observe that because $\theta^{k}_\ell(x)=\theta^{k}_0(x-\ell)$, it suffices to find $\partial_x \theta^{k}_0$ and then by shifting, the corresponding formulas for SOPWs with other shift indices follow easily. Now using formulas \eqref{equation:SOPWfourierLine1} and \eqref{equation:SOPWfourierLine2}, 
\begin{align}
\partial_x \theta^k_0 &= \frac{1}{\sqrt{L}}\sum_{\frac{(k-1)L}{2}<|n|<\frac{kL}{2}}(\sgn(n) i)^{k-1}\partial_x \phi_n+\frac{1}{\sqrt{2L}}\sum_{|n|=\frac{(k-1)L}{2},\frac{kL}{2}}(\sgn(n) i)^{k-1}\partial_x \phi_n \notag \\
& = \frac{1}{\sqrt{L}}\sum_{\frac{(k-1)L}{2}<|n|<\frac{kL}{2}}(\sgn(n) i)^{k-1} (\frac{i2\pi n}{L}) \phi_n+\frac{1}{\sqrt{2L}}\sum_{|n|=\frac{(k-1)L}{2},\frac{kL}{2}}(\sgn(n) i)^{k-1}  (\frac{i2\pi n}{L}) \phi_n. 
\label{equation:1derivativeSOPWk}\end{align}
Now from equations \eqref{equation:fourierSOPWLine2} and using lemma \ref{lemma:sumOmega},
\begin{align}
~ & \frac{1}{\sqrt{L}}\sum_{\frac{(k-1)L}{2}<|n|<\frac{kL}{2}}(\sgn(n) i)^{k-1} (\frac{i2\pi n}{L}) \phi_n \notag \\
= & \frac{1}{\sqrt{L}}\sum_{\frac{(k-1)L}{2}<|n|<\frac{kL}{2}}(\sgn(n) i)^{k-1} (\frac{i2\pi n}{L}) \frac{(- \sgn(n) i)^{k-1}}{\sqrt{L}}\sum_{j=0}^{L-1} \omega_j^n \theta_i^k \notag \\
= & \frac{i2\pi}{L^2}\sum_{j=0}^{L-1}\left( \sum_{\frac{(k-1)L}{2}<|n|<\frac{kL}{2}}n\omega_j^n \right) \theta_j^k \notag \\
= & \frac{\pi}{L} \sum_{j=0}^{L-1} a(j) \theta_j^k, \label{equation:1derivativePart1}
\end{align}
where 
\[  a(j)= \begin{cases}
0 \qquad &\hbox{if } j=0, \\
(-1)^{k}(2k-1) \cot(\pi j/L)   &\hbox{if $j$ is odd}, \\ 
\cot(\pi j/L)  &\hbox{otherwise}.
\end{cases} \]
On the other hand, equation \eqref{equation:fourierSOPWLine3} implies that for $|n|=\frac{kL}{2}$,
\begin{equation}  \phi_n(x)=\frac{(-\sgn(n) i)^{k-1}}{\sqrt{2L}}\left (\sum_{j=0}^{L-1}(-1)^{kj} \theta^k_j(x)- \sgn(n) i \sum_{j=0}^{L-1}(-1)^{kj} \theta^{k+1}_j(x) \right). \label{equation:phi=nk/2}\end{equation}
Therefore,
\begin{align}
~&\frac{1}{\sqrt{2L}}\sum_{|n|=\frac{kL}{2}} (\sgn(n) i)^{k-1}  (\frac{i2\pi n}{L}) \phi_n \notag \\
=& \frac{1}{2L}(\frac{i2\pi}{L})\sum_{|n|=\frac{kL}{2}}\left( n\sum_{j=0}^{L-1}(-1)^{kj} \theta^{k}_j- |n|i  \sum_{j=0}^{L-1}(-1)^{kj} \theta^{k+1}_j  \right) \notag\\
= & \frac{\pi}{L} k \sum_{j=0}^{L-1}(-1)^{kj}\theta^{k+1}_j. 
\label{equation:1derivativePart2}\end{align}
Again, equation \eqref{equation:fourierSOPWLine3} implies that for $|n|=\frac{(k-1)L}{2}$,
\begin{equation}  \phi_n(x)=\frac{(-\sgn(n) i)^{k-2}}{\sqrt{2L}}\left (\sum_{j=0}^{L-1}(-1)^{(k-1)j} \theta^{k-1}_j(x)- \sgn(n) i \sum_{j=0}^{L-1}(-1)^{(k-1)j} \theta^{k}_j(x) \right). \label{equation:phi=nk1/2} \end{equation}
Therefore,
\begin{align}
~&\frac{1}{\sqrt{2L}}\sum_{|n|=\frac{(k-1)L}{2}} (\sgn(n) i)^{k-1}  (\frac{i2\pi n}{L}) \phi_n \notag \\
=& \frac{1}{2L}(\frac{i2\pi}{L})\sum_{|n|=\frac{(k-1)L}{2}}\left( |n| i \sum_{j=0}^{L-1}(-1)^{(k-1)j} \theta^{k-1}_j+n  \sum_{j=0}^{L-1}(-1)^{(k-1)j} \theta^{k}_j  \right) \notag\\
= & -\frac{\pi}{L} (k-1) \sum_{j=0}^{L-1}(-1)^{(k-1)j}\theta^{k-1}_j. 
\label{equation:1derivativePart3}\end{align}
Substituting \eqref{equation:1derivativePart1}, \eqref{equation:1derivativePart2} and \eqref{equation:1derivativePart3} into equation \eqref{equation:1derivativeSOPWk} yields that
\[ \partial_x \theta^{k}_0= \frac{\pi}{L}\left[ -(k-1)\sum_{j=0}^{L-1} (-1)^{(k-1)j}\theta^{k-1}_j+\sum_{j=0}^{L-1} a(j) \theta^{k}_j+k\sum_{j=0}^{L-1} (-1)^{kj}\theta^{k+1}_j \right]. \]
This completest the proof. \qed

\textbf{Proof of theorem \ref{theorem:secondDerivative}:} Again observe that because $\theta^{k}_\ell(x)=\theta^{k}_0(x-\ell)$, it suffices to find $\partial_{xx} \theta^{k}_0$ and then by shifting, the corresponding formulas for SOPWs with other shift indices follow easily.  Now using formulas \eqref{equation:SOPWfourierLine1} and \eqref{equation:SOPWfourierLine2}, 
\begin{align}
~&\partial_{xx} \theta^k_0 \notag\\
=& \frac{1}{\sqrt{L}}\sum_{\frac{(k-1)L}{2}<|n|<\frac{kL}{2}}(\sgn(n) i)^{k-1}\partial_{xx} \phi_n+\frac{1}{\sqrt{2L}}\sum_{|n|=\frac{(k-1)L}{2},\frac{kL}{2}}(\sgn(n) i)^{k-1}\partial_{xx} \phi_n \notag \\
=& \frac{1}{\sqrt{L}}\sum_{\frac{(k-1)L}{2}<|n|<\frac{kL}{2}}(\sgn(n) i)^{k-1} (\frac{-4\pi^2 n^2}{L^2}) \phi_n+\frac{1}{\sqrt{2L}}\sum_{|n|=\frac{(k-1)L}{2},\frac{kL}{2}}(\sgn(n) i)^{k-1}  (\frac{-4\pi^2 n^2}{L^2}) \phi_n. 
\label{equation:2derivativeSOPWk}\end{align}
Now from equations \eqref{equation:fourierSOPWLine2} and using lemma \ref{lemma:sumOmega},
\begin{align}
~ & \frac{1}{\sqrt{L}}\sum_{\frac{(k-1)L}{2}<|n|<\frac{kL}{2}}(\sgn(n) i)^{k-1}  (\frac{-4\pi^2 n^2}{L^2})  \phi_n \notag \\
= & \frac{1}{\sqrt{L}}\sum_{\frac{(k-1)L}{2}<|n|<\frac{kL}{2}}(\sgn(n) i)^{k-1}  (\frac{-4\pi^2 n^2}{L^2})  \frac{(- \sgn(n) i)^{k-1}}{\sqrt{L}}\sum_{j=0}^{L-1} \omega_j^n \theta_i^k \notag \\
= & \frac{-4\pi^2}{L^3}\sum_{j=0}^{L-1}\left( \sum_{\frac{(k-1)L}{2}<|n|<\frac{kL}{2}}n^2\omega_j^n \right) \theta_j^k \notag \\
= & \frac{-\pi^2}{L^2} \sum_{j=0}^{L-1} \tilde{b}(j) \theta_j^k, \label{equation:2derivativePart1}
\end{align}
where 
\[  \tilde{b}(j)= \begin{cases}
(L-2)[(3k^2-3k+1)L-1]/3 \qquad &\hbox{if } j=0, \\
(-1)^{k}(4k-2)\csc^2(\pi j/L)-(-1)^k(2k-1)L   &\hbox{if $j$ is odd}, \\ 
2\csc^2(\pi j/L)-(k^2+(k-1)^2)L  &\hbox{otherwise}.
\end{cases} \]
On the other hand, from \eqref{equation:phi=nk/2},
\begin{align}
~&\frac{1}{\sqrt{2L}}\sum_{|n|=\frac{kL}{2}} (\sgn(n) i)^{k-1}  (\frac{-4\pi^2 n^2}{L^2}) \phi_n \notag \\
=& \frac{1}{2L}(\frac{-4\pi^2}{L^2})\sum_{|n|=\frac{kL}{2}}\left( n^2\sum_{j=0}^{L-1}(-1)^{kj} \theta^{k}_j- n^2 \sgn(n) i  \sum_{j=0}^{L-1}(-1)^{kj} \theta^{k+1}_j  \right) \notag\\
= & \frac{-\pi^2}{L^2} k^2L \sum_{j=0}^{L-1}(-1)^{kj}\theta^{k}_j. 
\label{equation:2derivativePart2}\end{align}
Also from \eqref{equation:phi=nk1/2},
\begin{align}
~&\frac{1}{\sqrt{2L}}\sum_{|n|=\frac{(k-1)L}{2}} (\sgn(n) i)^{k-1}   (\frac{-4\pi^2 n^2}{L^2}) \phi_n \notag \\
=& \frac{1}{2L} (\frac{-4\pi^2}{L^2})\sum_{|n|=\frac{(k-1)L}{2}}\left( n^2\sgn(n) i\sum_{j=0}^{L-1}(-1)^{(k-1)j} \theta^{k-1}_j+n^2  \sum_{j=0}^{L-1}(-1)^{(k-1)j} \theta^{k}_j  \right) \notag\\
= & \frac{-\pi^2}{L^2} (k-1)^2L \sum_{j=0}^{L-1}(-1)^{(k-1)j}\theta^{k}_j.
\label{equation:2derivativePart3}\end{align}
Substituting \eqref{equation:2derivativePart1}, \eqref{equation:2derivativePart2} and \eqref{equation:2derivativePart3} into equation \eqref{equation:2derivativeSOPWk} and simplifying yields that
\[ \partial_{xx} \theta^{k}_0= \frac{-\pi^2}{L^2} \sum_{j=0}^{L-1} b(j) \theta^{k}_j. \]
This completest the proof. \qed

\end{document}